\documentclass[english,10pt]{article}
\usepackage{amsmath}
\usepackage{graphicx}
\usepackage{amssymb}
\usepackage[usenames, dvipsnames]{color}
\usepackage{mathrsfs}
\usepackage{lineno}
\usepackage{amsfonts,epsfig}
\usepackage{amsfonts,amssymb,bm}
\usepackage{txfonts}
\usepackage{pifont}
\usepackage{tikz}
\usepackage[amsmath,thref,thmmarks]{ntheorem}

\usepackage{lineno}

\def\0{\emptyset}

\newtheorem{Theorem}{Theorem}[section]

\newtheorem{Example}{Example}[section]

\newtheorem{Proposition}{Proposition}[section]
\newtheorem{Remark}{Remark}[section]

\definecolor{blue}{rgb}{0.0,0.0,1.0}

\font\euler=eusm10
\def \M{\mbox{\euler M}}

\begin{document}

\title{High-order Uncertain Differential Equation and Its Application to Nuclear Reactors}

\author{
Hao Li\thanks{School of Mathematics, Renmin University of China, Beijing 100872, China. Email: hlimath@ruc.edu.cn,  wangyuqian@ruc.edu.cn},
Yuqian Wang$^*$ }
\date{}
\maketitle

\begin{abstract} High-order uncertain differential equation (HUDE) was introduced in literature. 
But the present method to solve a
HUDE is incorrect. In this paper, we will rigorously prove some comparion theorems of high-order differential equations, and present a 
method to solve a family of HUDE, including parameter estimation and hypothesis test. Then an application to nuclear reactor kinetics is given 
to illustrate the method.
\end{abstract}

\\

\vskip 0.35cm

\section{Introduction}
Uncertainty theory, founded by \cite{Liu2007}, has been developed into an axiomatic mathematical theory. 
Among the many theoretical branches of uncertainty theory, uncertain 
statistics has been the most popular and cutting-edge theoretical branch for the last few years. There are three 
important methods in uncertain statistics: uncertain time series analysis, uncertain regression analysis and 
uncertain differential equations. They have drawn attention of researchers with different background. Many achievements and 
results have been made in theory and applications. 

Uncertain differential equation was first introduced by \cite{Liu2008}. For those uncertain differential equations with no analytical solutions, 
\cite{Yao-Chen2013} presented a useful formula to calculate the inverse uncertainty distribution of the solution to an uncertain differential equation in terms of $\alpha$-path. This formula is known as "Yao-Chen Formula".  
\cite{Zhu2023} discussed the uncertain partial differential equations.
\cite{Ye2023} discussed the partial derivatives of uncertain field, and gave the integral form of an uncertain partial differential equations.

Uncertain differential equations has been widely applied to many fields, such as chemical reaction (\cite{TangandYang2021}), pharmacokinetics (\cite{LiuandYang2021}),
epidemic spread (\cite{Lio2021}), gas price (\cite{Mehrdoust2023}), China's population(\cite{Yang and Liu2023}), China's birth rates(\cite{YeandZheng2023}). Recently, some practical analysis in finance verify that, compared with 
stochastic differential equations, uncertain  differential 
equations are more suitable to fit the data, for example, Alibaba stock price (\cite{LiuandLiu2022}), currency exchange rate (\cite{LiuandYe2023}), 
interest rate (\cite{YangandKe2023}). 
When applying uncertain differential equations to practical problems, there are two core problems: how to estimate unknown
parameters in an uncertain differential equation based on the observed data, and how to test the fitness of an uncertain differential equation. 
\cite{LiuandLiu2022} introduced the concept of residuals of uncertain differential equations, and developed the method of moments estimation.
Soon, \cite{LiuandLiu2023-1} explored a modified maximum likelihood estimation, and \cite{LiuandLiu2023-2} presented the least squares estimation.
For evaluating an uncertain differential equation's goodness of fit, \cite{LiuandYe2023} intorduced the hypothesis test.  

For complex dynamic systems, such as spring vibration, pendulum swing, RLC circuit, nuclear reactor kinetics, high-order uncertain differential equations are required to characterize them.
\cite{Yao2016} initially proposed the HUDE. 
However, in the process of research, we found that there is a mistake in \cite{Yao2016} that the theorem of $\alpha$-path of HUDE is not strictly proved. Therefore, the mathod of solving a high-order uncertian differentail equation in \cite{Yao2016} is wrong.

In this paper, we focus on high-order uncertain differential equations. We will prove some comparion theorems of high-order ordinary differential equation, and rigorously prove a theorem of $\alpha$-path of HUDE. Then we present a method to solve a HUDE, including parameter estimation and hypothesis test.
As an application, the nuclear reactor kinetics under uncertain circumstance is discussed.
The remainder of this paper is organized as follows. In Section \ref{Preliminary}, we will prove comparison theorems and give the solution of  
a HUDE. In Section \ref{residuals}, we will introduce the concept of residual of a HUDE. 
We will discuss parameter estimation and hypothesis test in Section \ref{PE} and Section \ref{UHT}, respectively. In Section \ref{nuclear reactor},
the model of HUDE will be applied to nuclear reactor kinetics. Finally, some conclusions will be made in Section 7.


\section{Solution of a High-order Uncertain Differential Equation}\label{Preliminary}
An uncertain differential equation 
\begin{equation}\label{high}
\frac{\mathrm{d}^n X_t}{\mathrm{d}t^n}=f\bigg(t,X_t,\cdots,\frac{\mathrm{d}^{n-1} X_t}{\mathrm{d}t^{n-1}}\bigg)+\sum_{i=1}^{m}g_i\bigg(t,X_t,\cdots,\frac{\mathrm{d}^{n-1} X_t}{\mathrm{d}t^{n-1}}\bigg)\frac{\mathrm{d}C_{it}}{\mathrm{d}t}
\end{equation}
is called a high-order uncertain differential equation, where $f$ and $g_i\ (i=1,2,\cdots,m)$ are continuous functions, and $C_{1t},C_{2t},\cdots,C_{mt}$ are independent Liu processes.
Since most of the high-order uncertain differential equations cannot be solved analytically, the uncertainty distribution of the solution cannot be determined. So it is significantly important to figure out the inverse uncertainty distribution of the solution.

For a high-order uncertain differential equation \eqref{high}, let $X_t^{\alpha}$ be the solution of the corresponding ordinary differential equation:
	\begin{equation}\label{path}
	\frac{\mathrm{d}^n X_t^{\alpha}}{\mathrm{d}t^n}=f\bigg(t,X_t^{\alpha},\cdots,\frac{\mathrm{d}^{n-1} X_t^{\alpha}}{\mathrm{d}t^{n-1}}\bigg)+\sum_{i=1}^{m}\bigg|g_i\bigg(t,X_t^{\alpha},\cdots,\frac{\mathrm{d}^{n-1} X_t^{\alpha}}{\mathrm{d}t^{n-1}}\bigg)\bigg|\Phi^{-1}(\alpha)
	\end{equation}
	where
	$$
	\Phi^{-1}(\alpha)=\frac{\sqrt{3}}{\pi}\mathrm{ln}\frac{\alpha}{1-\alpha},\qquad 0<\alpha<1.
	$$

In \cite{Yao2016}, $X_t^{\alpha}$ is proved to be the $\alpha$-path of $X_t$, where $X_t$ is the solution of \eqref{high}.
The relevant conclusion is listed below. \\

\noindent {\bf Conclusion:}(\cite{Yao2016}, Page 144, Theorem 10.1)
The solution $X_t$ of a high-order uncertain differential equation
\begin{equation*}
\frac{\mathrm{d}^n X_t}{\mathrm{d}t^n}=f\bigg(t,X_t,\cdots,\frac{\mathrm{d}^{n-1} X_t}{\mathrm{d}t^{n-1}}\bigg)+g\bigg(t,X_t,\cdots,\frac{\mathrm{d}^{n-1} X_t}{\mathrm{d}t^{n-1}}\bigg)\frac{\mathrm{d}C_{t}}{\mathrm{d}t}
\end{equation*}
is a contour process with an $\alpha$-path $X_t^{\alpha}$ that solves the corresponding high-order ordinary differential equation
\begin{equation*}
\frac{\mathrm{d}^n X_t^{\alpha}}{\mathrm{d}t^n}=f\bigg(t,X_t^{\alpha},\cdots,\frac{\mathrm{d}^{n-1} X_t^{\alpha}}{\mathrm{d}t^{n-1}}\bigg)+\bigg|g\bigg(t,X_t^{\alpha},\cdots,\frac{\mathrm{d}^{n-1} X_t^{\alpha}}{\mathrm{d}t^{n-1}}\bigg)\bigg|\Phi^{-1}(\alpha)
\end{equation*}
where
$$
\Phi^{-1}(\alpha)=\frac{\sqrt{3}}{\pi}\mathrm{ln}\frac{\alpha}{1-\alpha},\qquad 0<\alpha<1
$$	
is the inverse uncertainty distribution of standard normal uncertain variables. In other words, 
$$\M\big\{X_t\le X_t^{\alpha}, \forall t\big\}=\alpha,$$
$$\M\big\{X_t\ge X_t^{\alpha}, \forall t\big\}=1-\alpha.$$ 

However, there is a drawback in its proof. Here is a counterexample.

\begin{Example}
	Consider the following uncertain differential equation when $t\ge 0$
	\begin{equation}\label{ude-1}
	\left\{
	\begin{array}{ll}
	\vspace{2mm}
    &\displaystyle\frac{\mathrm{d}^2 X_t}{\mathrm{d}t^2}=-X_t+e^{-t}\frac{\mathrm{d}C_{t}}{\mathrm{d}t}\\
    \vspace{2mm}
    &\displaystyle X_{0}=a\\
    &\displaystyle\frac{\mathrm{d}X_{t}}{\mathrm{d}t}\bigg|_{t=0}=b
	\end{array}
	\right.
	\hspace{30pt}
	\end{equation} 
	where $a$ and $b$ are constants.
	Suppose that $\Psi_t^{-1}(\alpha)$ is the inverse uncertainty distribution of $X_t$.
    If the conclusion is correct, 
    then we have
    \begin{equation}\label{Psi}
    \Psi_t^{-1}(\alpha)=X_t^{\alpha}.
    \end{equation}
    Note that $\Psi_t^{-1}(\alpha)$ is an increasing function with respect to $\alpha$. 
    Choose $\alpha_1, \alpha_2$ such that $0<\alpha_1<\alpha_2<1$.
    Let $p_1=\Phi^{-1}(\alpha_1)$, $p_2=\Phi^{-1}(\alpha_2)$, and then $p_1<p_2$.
    Assume $X_t^{\alpha_i}(i=1,2)$ is the solution of 
    \begin{equation}\label{ude-2}
    \left\{
    \begin{array}{ll}
    \vspace{2mm}
    &\displaystyle\frac{\mathrm{d}^2 X_t}{\mathrm{d}t^2}=-X_t+p_ie^{-t}\\
    \vspace{2mm}
    &\displaystyle X_{0}=a\\
    &\displaystyle\frac{\mathrm{d}X_{t}}{\mathrm{d}t}\bigg|_{t=0}=b.
    \end{array}
    \right.
    \hspace{30pt}
    \end{equation}
    It is easy to check that
    $$
    X_t^{\alpha_i}=\Big(a-\frac{p_i}{2}\Big)\cos{t}+\Big(b+\frac{p_i}{2}\Big)\sin{t}+\frac{p_i}{2}e^{-t},\quad i=1,2.
    $$
    The graphs of $X_t^{\alpha}$ of some different values of $\alpha$ are shown in Fig. \ref{graph6}.
    Then
    \begin{equation}\label{eq:compare2}
    X_t^{\alpha_1}-X_t^{\alpha_2}=\frac{1}{2}(p_1-p_2)\bigg(\sqrt{2}\sin\Big(t-\frac{\pi}{4}\Big)+e^{-t}\bigg).
    \end{equation}

    \begin{figure}[htbp]
    	\centering
    	\includegraphics[scale=0.5]{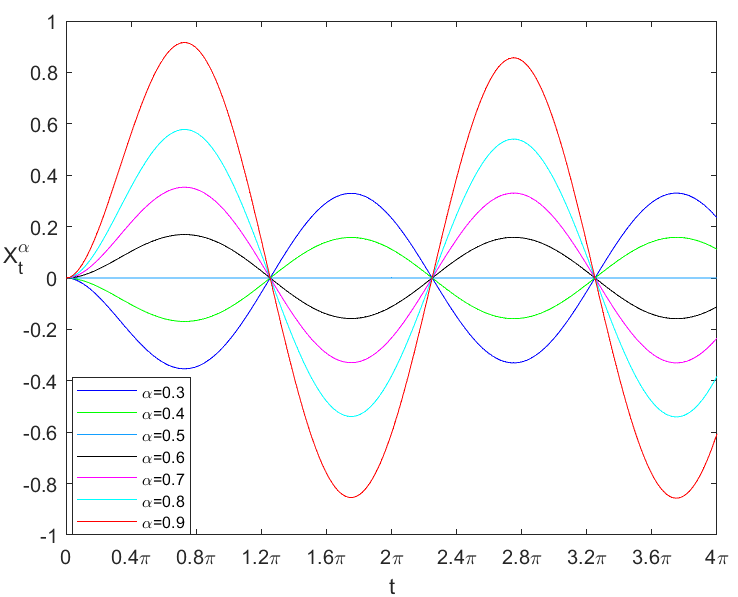}
    	\caption{The graphs of $X_t^{\alpha}$ of some different values of $\alpha$}
    	\label{graph6}
    \end{figure}
    
    So
    $$
    X_t^{\alpha_1}-X_t^{\alpha_2}>0
    $$
    for 
    $$t\in \bigg(\displaystyle-\frac{3\pi}{4}+2k\pi, \frac{\pi}{4}+2k\pi\bigg)\ (k\in \mathbb{N}^{+}).$$ 
    Thus, by Eq.\eqref{Psi}, 
    $$
    \Psi_t^{-1}(\alpha_1)=X_t^{\alpha_1}>X_t^{\alpha_2}=\Psi_t^{-1}(\alpha_2).
    $$
    This is a contradiction to the fact that $\Psi_t^{-1}(\alpha)$ is an increasing function.
    Therefore, this conclusion is incorrect.
\end{Example}

Next, we will prove two comparison theorems of ordinary differential equation, and give a sufficient condition for $X_t^{\alpha}$ to be the $\alpha$-path of $X_t$.

\begin{Theorem}\label{compare1}
	 Let $f\big(t,z_0,z_1,\cdots,z_{n-1}\big)$ and $g\big(t,z_0,z_1,\cdots,z_{n-1}\big)$ be two functions on $D=[t_0, a]\!\times\! \mathbb{D}_0\subset\mathbb{R}^{n\!+\!1}$ (where $\mathbb{D}_0\  \mbox{is an n-dimensional bounded region})$ satisfying local Lipschitz conditions in $z_0,z_1,\cdots,z_{n-1}$. Assume that 
	 $f\big(t,z_0,z_1,$ $\cdots,z_{n-1}\big) <g\big(t,z_0,z_1,\cdots,z_{n-1}\big)$ on $D$, and that at least one of $f\big(t,z_0,z_1,\cdots,z_{n-1}\big)$ and $g\big(t,z_0,z_1,\cdots,z_{n-1}\big)$ is a monotonically increasing function with respect to $z_0,z_1,\cdots,$ $z_{n-2}$.
	 If $\psi(t)$ and $\Psi(t)$ are solutions of
	\begin{equation*}
	\left\{
	\begin{aligned}
	&\frac{\mathrm{d}^{n} X_t}{\mathrm{d}t^{n}}=f\bigg(t,X_t,\cdots,\frac{\mathrm{d}^{n\!-\!1} X_t}{\mathrm{d}t^{n\!-\!1}}\bigg)\\
	&X(t_0)=y_0, \frac{\mathrm{d}X_{t}}{\mathrm{d}t}\bigg|_{t=t_0}=y'_0, \cdots, \frac{\mathrm{d}^{n-1} X_t}{\mathrm{d}t^{n-1}}\bigg|_{t=t_0}=y_0^{(n-1)}
	\end{aligned}
	\right.
	\end{equation*}
	and 
	$$
	\left\{
	\begin{aligned}
	&\frac{\mathrm{d}^{n} X_t}{\mathrm{d}t^{n}}=g\bigg(t,X_t,\cdots,\frac{\mathrm{d}^{n\!-\!1} X_t}{\mathrm{d}t^{n\!-\!1}}\bigg)\\
	&X(t_0)=y_0, \frac{\mathrm{d}X_{t}}{\mathrm{d}t}\bigg|_{t=t_0}=y'_0, \cdots, \frac{\mathrm{d}^{n-1} X_t}{\mathrm{d}t^{n-1}}\bigg|_{t=t_0}=y_0^{(n-1)},
	\end{aligned}
	\right.
	$$
	respectively, then 
	\begin{equation*}	
	\psi(t)<\Psi(t),\ \forall t\in [t_0, a].
	\end{equation*}
\end{Theorem}	

\noindent$\mathbf{Proof.}$ Let
$$
\phi(t)=\Psi(t)-\psi(t), \quad t\in [t_0, a].
$$
Then $\phi(t_0)=\phi'(t_0)=\cdots=\phi^{(n-1)}(t_0)=0,$ and $ \phi^{(n)}(t_0)>0$. 
By the sign-preserving property of derivatives, there exists $\delta >t_0$ such that 
\begin{equation}\label{eq-sigma}
\phi(t)>0,\ \phi'(t)>0,\cdots, \phi^{(n-1)}(t)>0,\ \forall t\in (t_0, \delta).
\end{equation} 

By way of contradiction, suppose that $\phi(t)$ is not always positive when $t\in [t_0, a]$.
Let 
$$
t_1=\min\big\{t~|~\phi(t)=0,\ t\in [t_0, a]\big\}.
$$
Then 
\begin{equation}
\delta \le t_1, \mbox{ and }\phi(t)>0,\ \forall t\in(t_0, t_1).
\end{equation}
By the choice of $t_1$ and the above inequalities, 
\begin{equation}\label{eq:t1}
\phi'(t_1)\le 0.
\end{equation}

\begin{figure}[htbp]
	\centering
	\includegraphics[scale=0.5]{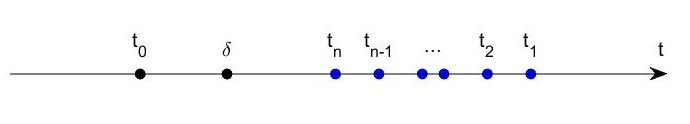}
	\caption{$\{t_1, t_2, \cdots, t_n\}$ is 
decreasing with a lower bound $\delta$}
	\label{graph0}
\end{figure}

Next, we define $t_2, t_3, \cdots, t_{n-1}$ and $t_n$ one by one in the following way.
For $i\ (2\le i\le n)$, suppose $t_{i-1}$ is defined such that 
$$
\delta \le t_{i-1},\quad \phi^{(i-1)}(t_{i-1})\le 0
$$
and 
$$
\phi(t)>0,\phi'(t)>0,\cdots,\phi^{(i-2)}(t)>0,\ \forall t\in(t_0, t_{i-1}). 
$$
Let
\begin{equation}\label{def ti}
t_i=\min\Big\{t~\big|~\phi^{(i-1)}(t)=0,\ t\in [t_0, a]\Big\}.
\end{equation}
Then $\delta\le t_i\le t_{i-1}$ and
\begin{equation}\label{ti}
\phi(t)>0,\phi'(t)>0,\cdots,\phi^{(i-1)}(t)>0,\ \forall t\in(t_0, t_i).
\end{equation}
By the choice of $t_i$ and the above inequalities,
\begin{equation}\label{eq:tn}
\phi^{(i)}(t_i)\le 0.
\end{equation}
Repeat this process until we have $t_1, t_2,\cdots, t_n$.
By the definition of $t_i$ ($1\le i \le n$), the sequence $\{t_1, t_2, \cdots, t_n\}$ is 
decreasing with a lower bound $\delta$ (Fig. \ref{graph0}).

Let $k$ be the minimum number such that $t_k=t_n$.
Then $t_{k-1}>t_k=t_{k+1}=\cdots =t_n$. By the choice of $t_n$,
$$
\phi(t)>0,\phi'(t)>0,\cdots,\phi^{(n-1)}(t)>0,\ \forall t\in(t_0, t_k)
$$ 
and
\begin{equation}\label{tk}
\phi^{(n)}(t_k)\le 0.
\end{equation}

As $t_k=t_{k+1}=\dots =t_n$, $\phi^{(i)}(t_k)=0$ for $k-1\le i \le n-1$ by Eq.(\ref{def ti}).
That is,
\begin{equation}\label{eva1}
\Psi^{(i)}(t_k)=\psi^{(i)}(t_k), \mbox{ for } k-1\le i \le n-1.
\end{equation}
As $t_k<t_{k-1}$, by Eq.\eqref{ti}, $\phi^{(i)}(t_k)>0$ for $0\le i \le k-2$, i.e.
\begin{equation}\label{eva2}
\Psi^{(i)}(t_k)>\psi^{(i)}(t_k), \mbox{ for } 0\le i \le k-2.
\end{equation}
If $g\big(t,z_0,z_1,\cdots,z_{n-1}\big)$ is a monotonically increasing function with respect to $z_0,\cdots,$ $z_{n-2}$, by Eq. \eqref{eva1} and \eqref{eva2},
$$\begin{aligned}
\phi^{(n)}(t_k)&=\Psi^{(n)}(t_k)-\psi^{(n)}(t_k)\\
&=g\Big(t_k,\Psi(t_k),\cdots,\Psi^{(n-1)}(t_k)\Big)-f\Big(t_k,\psi(t_k),\cdots,\psi^{(n-1)}(t_k)\Big)\\
&\ge g\Big(t_k,\psi(t_k),\cdots,\psi^{(n-1)}(t_k)\Big)-f\Big(t_k,\psi(t_k),\cdots,\psi^{(n-1)}(t_k)\Big)\\
&>0.
\end{aligned}
$$
If $f\big(t,z_0,z_1,\cdots,z_{n-1}\big)$ is a monotonically increasing function with respect to $z_0,\cdots,$ $z_{n-2}$, by Eq. \eqref{eva1} and \eqref{eva2},
$$\begin{aligned}
\phi^{(n)}(t_k)&=g\Big(t_k,\Psi(t_k),\cdots,\Psi^{(n-1)}(t_k)\Big)-f\Big(t_k,\psi(t_k),\cdots,\psi^{(n-1)}(t_k)\Big)\\
&\ge g\Big(t_k,\Psi(t_k),\cdots,\Psi^{(n-1)}(t_k)\Big)-f\Big(t_k,\Psi(t_k),\cdots,\Psi^{(n-1)}(t_k)\Big)\\
&>0.
\end{aligned}
$$
This contradicts to Eq.\eqref{tk}.
Therefore  
this assumption is not valid and the conclusion $\phi(t)>0$ for $ t\in [t_0, a]$  is true, i.e.
$$
\psi(t)<\Psi(t),\ \forall t\in [t_0, a].
$$
The theorem is proved.

\begin{Theorem}\label{compare2}
	Let $f\big(t,z_0,z_1,\cdots,z_{n-1}\big)$ and $g\big(t,z_0,z_1,\cdots,z_{n-1}\big)$ be two functions on $D=[t_0, a]\!\times\! \mathbb{D}_0\subset\mathbb{R}^{n\!+\!1}$ (where $\mathbb{D}_0\  \mbox{is an n-dimensional bounded region})$ satisfying local Lipschitz conditions in $z_0,z_1,\cdots,z_{n-1}$. Assume that 
	$f\big(t,z_0,z_1,$ $\cdots,z_{n-1}\big) \le g\big(t,z_0,z_1,\cdots,z_{n-1}\big)$ on $D$, and that at least one of $f\big(t,z_0,z_1,\cdots,z_{n-1}\big)$ and $g\big(t,z_0,z_1,\cdots,z_{n-1}\big)$ is a monotonically increasing function with respect to $z_0,z_1,\cdots,$ $z_{n-2}$.
	If $\psi(t)$ and $\Psi(t)$ are solutions of
	\begin{equation*}
	\left\{
	\begin{aligned}
	&\frac{\mathrm{d}^{n} X_t}{\mathrm{d}t^{n}}=f\bigg(t,X_t,\cdots,\frac{\mathrm{d}^{n\!-\!1} X_t}{\mathrm{d}t^{n\!-\!1}}\bigg)\\
	&X(t_0)=y_0, \frac{\mathrm{d}X_{t}}{\mathrm{d}t}\bigg|_{t=t_0}=y'_0, \cdots, \frac{\mathrm{d}^{n-1} X_t}{\mathrm{d}t^{n-1}}\bigg|_{t=t_0}=y_0^{(n-1)}
	\end{aligned}
	\right.
	\end{equation*}
	and 
	$$
	\left\{
	\begin{aligned}
	&\frac{\mathrm{d}^{n} X_t}{\mathrm{d}t^{n}}=g\bigg(t,X_t,\cdots,\frac{\mathrm{d}^{n\!-\!1} X_t}{\mathrm{d}t^{n\!-\!1}}\bigg)\\
    &X(t_0)=y_0, \frac{\mathrm{d}X_{t}}{\mathrm{d}t}\bigg|_{t=t_0}=y'_0, \cdots, \frac{\mathrm{d}^{n-1} X_t}{\mathrm{d}t^{n-1}}\bigg|_{t=t_0}=y_0^{(n-1)},
	\end{aligned}
	\right.
	$$
	respectively, then 
	\begin{equation*}	
	\psi(t)\le\Psi(t),\ \forall t\in [t_0, a].
	\end{equation*}
\end{Theorem}
\noindent$\mathbf{Proof.}$
Let $\{\varepsilon_i\}\ (i=1,2,\cdots)$ be a monotonically decreasing sequence of positive numbers such that $\lim\limits_{i\rightarrow\infty}\varepsilon_i=0$. 
Due to
$$f\big(t,z_0,z_1,\cdots,z_{n-1}\big) \le g\big(t,z_0,z_1,\cdots,z_{n-1}\big)$$
for $\big(t,z_0,z_1,\cdots,z_{n-1}\big)\in D$, we have
$$f\big(t,z_0,z_1,\cdots,z_{n-1}\big)-\varepsilon_i < g\big(t,z_0,z_1,\cdots,z_{n-1}\big)$$
for $\big(t,z_0,z_1,\cdots,z_{n-1}\big)\in D$, where $i=1,2,\cdots$.

Consider equations 
\begin{equation}\label{i}
\left\{
\begin{aligned}
&\frac{\mathrm{d}^{n} X_t}{\mathrm{d}t^{n}}=f\bigg(t,X_t,\cdots,\frac{\mathrm{d}^{n\!-\!1} X_t}{\mathrm{d}t^{n\!-\!1}}\bigg)-\varepsilon_i\\
&X(t_0)=y_0, \frac{\mathrm{d}X_{t}}{\mathrm{d}t}\bigg|_{t=t_0}=y'_0, \cdots, \frac{\mathrm{d}^{n-1} X_t}{\mathrm{d}t^{n-1}}\bigg|_{t=t_0}=y_0^{(n-1)}.
\end{aligned}
\right.
\end{equation}
From the existence and uniqueness theorem of solution, the initial value problem \eqref{i} has exactly one solution $\psi_i(t)\ (i=1,2,\cdots)$ in the interval $t_0\le t\le a$.
According to Theorem\ \ref{compare1}, we can get
$$
\psi_i(t)<\Psi(t)\ (i=1,2,\cdots),\ \forall t\in [t_0, a].
$$
Note that 
$$
\begin{aligned}
\big| \psi_i(t_1)-\psi_i(t_2)\big|&\le\int_{t_1}^{t_2}\idotsint_{\mathbb{D}_0} \bigg|f\bigg(t,X_t,\cdots,\frac{\mathrm{d}^{n\!-\!1} X_t}{\mathrm{d}t^{n\!-\!1}}\bigg)-\varepsilon_i \bigg|\,dX_t \dots dt\\
&\le \int_{t_1}^{t_2}\bigg(\idotsint_{\mathbb{D}_0} M\,dX_t \dots\bigg) dt\\
&\le M\cdot \left\|\mathbb{D}_0\right\|\cdot \left|t_2-t_1\right|
\end{aligned}
$$
$\big(t_1, t_2\in [t_0, a], i\!=\!1,2,\!\cdots\!\big)$, where $$M\!=\!\max\limits_{\big(t,X_t,\cdots,\frac{\mathrm{d}^{n\!-\!1} X_t}{\mathrm{d}t^{n\!-\!1}}\big)\in D}\left|f\bigg(t,X_t,\!\cdots\!,\frac{\mathrm{d}^{n\!-\!1} X_t}{\mathrm{d}t^{n\!-\!1}}\bigg)\!-\!\varepsilon_i\right|.$$ Therefore, $\psi_i(t)\ (i=1,2,\cdots)$ is uniformly bounded and equally continuous on $[t_0, a]$ based on local Lipschitz condition.
It follows from Ascoli Lemma that $\psi_i(t)\ (i=1,2,\cdots)$ has uniformly convergent subsequences on region $D$.
Then
$$
\lim_{i\rightarrow\infty}\psi_i(t)=\psi(t).
$$
Consequently, 
$$
\psi(t)=\lim_{i\rightarrow\infty}\psi_i(t)\leq\Psi(t),\ \forall t\in [t_0, a].
$$
The theorem is proved.\\

Next, we will prove the theorem of $\alpha$-path for HUDE.

\begin{Theorem}\label{main Formula}
	Let $X_t$ and $X_t^\alpha$ be the solution of 
	\begin{equation*}
	\frac{\mathrm{d}^n X_t}{\mathrm{d}t^n}=f\bigg(t,X_t,\frac{\mathrm{d} X_t}{\mathrm{d}t}, \cdots,\frac{\mathrm{d}^{n-1} X_t}{\mathrm{d}t^{n-1}}\bigg)+\sum_{i=1}^{m}g_i\bigg(t,X_t,\frac{\mathrm{d} X_t}{\mathrm{d}t}, \cdots,\frac{\mathrm{d}^{n-1} X_t}{\mathrm{d}t^{n-1}}\bigg)\frac{\mathrm{d}C_{it}}{\mathrm{d}t}
	\end{equation*}
	and
	\begin{equation*}
    \frac{\mathrm{d}^n X_t^{\alpha}}{\mathrm{d}t^n}=f\bigg(t,X_t^{\alpha},\cdots,\frac{\mathrm{d}^{n-1} X_t^{\alpha}}{\mathrm{d}t^{n-1}}\bigg)+\sum_{i=1}^{m}\bigg|g_i\bigg(t,X_t^{\alpha},\cdots,\frac{\mathrm{d}^{n-1} X_t^{\alpha}}{\mathrm{d}t^{n-1}}\bigg)\bigg|\Phi^{-1}(\alpha)
    \end{equation*}
    respectively, where
    $$\Phi^{-1}(\alpha)=\frac{\sqrt{3}}{\pi}\mathrm{ln}\frac{\alpha}{1-\alpha},\qquad 0<\alpha<1.$$ 
	If $$f\bigg(t,X_t^{\alpha}, \frac{\mathrm{d} X_t^{\alpha}}{\mathrm{d}t}, \cdots,\frac{\mathrm{d}^{n-1} X_t^{\alpha}}{\mathrm{d}t^{n-1}}\bigg)+\sum\limits_{i=1}\limits^{m}\bigg|g_i\bigg(t,X_t^{\alpha}, \frac{\mathrm{d} X_t^{\alpha}}{\mathrm{d}t},\cdots,\frac{\mathrm{d}^{n-1} X_t^{\alpha}}{\mathrm{d}t^{n-1}}\bigg)\bigg|\Phi^{-1}(\alpha)$$
	is a monotonically increasing function with respect to $$X_t^{\alpha}, \frac{\mathrm{d} X_t^{\alpha}}{\mathrm{d}t}, \cdots,\frac{\mathrm{d}^{n-2} X_t^{\alpha}}{\mathrm{d}t^{n-2}},$$ then $X_t^{\alpha}$ is the $\alpha$-path of $X_t$, i.e., 
\begin{equation*}
\begin{aligned}
&\M\big\{X_t\le X_t^{\alpha},\forall t\big\}=\alpha,\\
&\M\big\{X_t>X_t^{\alpha},\forall t\big\}=1-\alpha.
\end{aligned}
\end{equation*}	
\end{Theorem}
\noindent $\mathbf{Proof.}$ Given $\alpha\in(0,1)$, for each $X_t^{\alpha}$, we can divide it into two parts,
$$
T_i^+=\bigg\{t\bigg|g_i\bigg(t,X_t^{\alpha},\cdots,\frac{\mathrm{d}^{n-1} X_t^{\alpha}}{\mathrm{d}t^{n-1}}\bigg)\ge0\bigg\}, 
$$
$$
T_i^-=\bigg\{t\bigg|g_i\bigg(t,X_t^{\alpha},\cdots,\frac{\mathrm{d}^{n-1} X_t^{\alpha}}{\mathrm{d}t^{n-1}}\bigg)<0\bigg\},
$$
$i=1,2,\cdots,m$. It is obvious that $T_i^+\cap T_i^-=\varnothing$ and $T_i^+\cup T_i^-=[0,+\infty)$ for each $i$, $1\le i \le m$.

Next, we define
$$
\begin{aligned}
&\Lambda_{i1}^+=\left\{\gamma \bigg|\frac{\mathrm{d}C_{it}(\gamma)}{\mathrm{d}t}\le \Phi^{-1}(\alpha)\ \mbox{for any}\ t\in T_i^+\right\},\\
&\Lambda_{i1}^-=\left\{\gamma \bigg|\frac{\mathrm{d}C_{it}(\gamma)}{\mathrm{d}t}\ge \Phi^{-1}(1-\alpha)\ \mbox{for any}\ t\in T_i^-\right\},
\end{aligned}
$$
$i=1,2,\cdots,m$, where 
$$
\Phi^{-1}(\alpha)=\frac{\sqrt{3}}{\pi}\mathrm{ln}\frac{\alpha}{1-\alpha}. 
$$
Because $T_i^+$ and $T_i^-$ are disjoint sets, and $C_{1t}, \cdots,C_{mt} $ are independent increment processes, we have
$$
\M\{\Lambda_{i1}^+\}=\alpha,\quad \M\{\Lambda_{i1}^-\}=\alpha,\quad \M\{\Lambda_{i1}^+\cap \Lambda_{i1}^-\}=\alpha,
$$
 $i=1,2,\cdots,m$.
For any $\gamma\in\Lambda_{i1}^+\cap \Lambda_{i1}^-$, it is apparent that for any $t$, 
$$
\begin{aligned}
&g_i\left(t,X_t(\gamma),\cdots,\frac{\mathrm{d}^{n-1} X_t(\gamma)}{\mathrm{d}t^{n-1}}\right)\frac{\mathrm{d}C_{it}(\gamma)}{\mathrm{d}t}\\
&\le 
\left|g_i\left(t,X_t^{\alpha},\cdots,\frac{\mathrm{d}^{n-1} X_t^{\alpha}}{\mathrm{d}t^{n-1}}\right)\right|\Phi^{-1}(\alpha),
\end{aligned}
$$
$i=1,2,\cdots,m$.
Let $\Lambda_1^+\cap\Lambda_1^-=\bigcap\limits_{i=1}^{m}(\Lambda_{i1}^+\cap \Lambda_{i1}^-)$, $i=1,2,\cdots,m$. 

Since $C_{1t},C_{2t},\cdots,C_{mt}$ are independent and $\M\{\Lambda_{i1}^+\cap \Lambda_{i1}^-\}=\alpha$, $i=1,2,\cdots,m$, we have
$$
\M\{\Lambda_1^+\cap\Lambda_1^-\}=\M\left\{\bigcap\limits_{i=1}^{m}(\Lambda_{i1}^+\cap \Lambda_{i1}^-)\right\}
=\bigwedge_{1\le i\le m}\M\big\{\Lambda_{i1}^+\cap \Lambda_{i1}^-\big\}=\alpha.
$$
So, for any $\gamma\in \Lambda_1^+\cap\Lambda_1^-$, we get for any $t$,
$$
\begin{aligned}
&\sum_{i=1}^{m}g_i\left(\!t,X_t(\!\gamma),\cdots,\frac{\mathrm{d}^{n\!-\!1} X_t(\!\gamma)}{\mathrm{d}t^{n\!-\!1}}\!\right)\frac{\mathrm{d}C_{it}(\gamma)}{\mathrm{d}t}\\
&\le \sum_{i=1}^{m}\left|g_i\left(\!t,X_t^{\alpha},\cdots,\frac{\mathrm{d}^{n\!-\!1} X_t^{\alpha}}{\mathrm{d}t^{n\!-\!1}}\!\right)\right|\Phi^{\!-\!1}(\alpha),
\end{aligned}
$$
i.e.
$$
\begin{aligned}
&f\left(t,X_t(\gamma),\cdots,\frac{\mathrm{d}^{n-1} X_t(\gamma)}{\mathrm{d}t^{n-1}}\right)+ \sum_{i=1}^{m}g_i\left(t,X_t(\gamma),\cdots,\frac{\mathrm{d}^{n-1} X_t(\gamma)}{\mathrm{d}t^{n-1}}\right)\frac{\mathrm{d}C_{it}(\gamma)}{\mathrm{d}t}\\
&\le f\left(t,X_t^{\alpha},\cdots,\frac{\mathrm{d}^{n-1} X_t^{\alpha}}{\mathrm{d}t^{n-1}}\right)+\sum_{i=1}^{m}\left|g_i\left(t,X_t^{\alpha},\cdots,\frac{\mathrm{d}^{n-1} X_t^{\alpha}}{\mathrm{d}t^{n-1}}\right)\right|\Phi^{-1}(\alpha).
\end{aligned}
$$
Since $$f\left(t,X_t^{\alpha},\cdots,\frac{\mathrm{d}^{n-1} X_t^{\alpha}}{\mathrm{d}t^{n-1}}\right)+\sum\limits_{i=1}\limits^{m}\left|g_i\left(t,X_t^{\alpha},\cdots,\frac{\mathrm{d}^{n-1} X_t^{\alpha}}{\mathrm{d}t^{n-1}}\right)\right|\Phi^{-1}(\alpha)$$
is a monotonically increasing function with respect to $$X_t^{\alpha}, \frac{\mathrm{d} X_t^{\alpha}}{\mathrm{d}t}, \cdots,\frac{\mathrm{d}^{n-2} X_t^{\alpha}}{\mathrm{d}t^{n-2}},$$ according to Theorem\ \ref{compare2}, we have
$$
X_t\le X_t^{\alpha}, \ \forall t.
$$
Note that $\Lambda_1^+\cap\Lambda_1^-\subset \{X_t\le X_t^{\alpha},\forall t\}$, we can get
\begin{equation}\label{first}
\M\{X_t\le X_t^{\alpha},\forall t\}\ge \M\{\Lambda_1^+\cap\Lambda_1^-\}=\alpha.
\end{equation}

Next, let 
$$
\begin{aligned}
&\Lambda_{i2}^+=\left\{\gamma \bigg|\frac{\mathrm{d}C_{it}(\gamma)}{\mathrm{d}t}> \Phi^{-1}(\alpha)\ \mbox{for any}\ t\in T_i^+\right\},\\
&\Lambda_{i2}^-=\left\{\gamma \bigg|\frac{\mathrm{d}C_{it}(\gamma)}{\mathrm{d}t}< \Phi^{-1}(1-\alpha)\ \mbox{for any}\ t\in T_i^-\right\},
\end{aligned}
$$
$i=1,2,\cdots,m$. Because $T_i^+$ and $T_i^-$ are disjoint sets and $C_{1t},\cdots,C_{mt}$ are independent increment processes, we get
$$
\M\{\Lambda_{i2}^+\}=1-\alpha,\quad \M\{\Lambda_{i2}^-\}=1-\alpha,\quad \M\{\Lambda_{i2}^+\cap \Lambda_{i2}^-\}=1-\alpha,
$$
$i=1,2,\cdots,m$. Considering $\forall\gamma\in\Lambda_{i2}^+\cap \Lambda_{i2}^-$, it is apparent that for any $t$,
$$
\begin{aligned}
&g_i\left(t,X_t(\gamma),\cdots,\frac{\mathrm{d}^{n-1} X_t(\gamma)}{\mathrm{d}t^{n-1}}\right)\frac{\mathrm{d}C_{it}(\gamma)}{\mathrm{d}t}\\
&> \left|g_i\left(t,X_t^{\alpha},\cdots,\frac{\mathrm{d}^{n-1} X_t^{\alpha}}{\mathrm{d}t^{n-1}}\right)\right|\Phi^{-1}(\alpha),
\end{aligned}
$$
$i=1,2,\cdots,m$.

Let $\Lambda_2^+\cap\Lambda_2^-=\bigcap\limits_{i=1}^{m}(\Lambda_{i2}^+\cap \Lambda_{i2}^-),i=1,2,\cdots,m$. 
Since $C_{1t},C_{2t},\cdots,C_{mt}$ are independent and $\M\{\Lambda_{i2}^+\cap \Lambda_{i2}^-\}\!=\!1-\alpha, i=1,2,\!\cdots\!,m$, we have
$$
\M\{\Lambda_2^+\cap\Lambda_2^-\}=\M\left\{\bigcap\limits_{i=1}^{m}(\Lambda_{i2}^+\cap \Lambda_{i2}^-)\right\}
=\bigwedge_{1\le i\le m}\M\big\{\Lambda_{i2}^+\cap \Lambda_{i2}^-\big\}=1-\alpha.
$$
So, for any $\gamma\in \Lambda_2^+\cap\Lambda_2^-$, we get for any $t$,
$$
\begin{aligned}
&\sum_{i=1}^{m}g_i\left(\!t,X_t(\!\gamma),\cdots,\frac{\mathrm{d}^{n\!-\!1} X_t(\!\gamma)}{\mathrm{d}t^{n\!-\!1}}\!\right)\frac{\mathrm{d}C_{it}(\gamma)}{\mathrm{d}t}\\
&> \sum_{i=1}^{m}\left|g_i\left(\!t,X_t^{\alpha},\cdots,\frac{\mathrm{d}^{n\!-\!1} X_t^{\alpha}}{\mathrm{d}t^{n\!-\!1}}\!\right)\right|\Phi^{\!-\!1}(\alpha).
\end{aligned}
$$
According to Theorem\ \ref{compare1}, we have
$$
X_t> X_t^{\alpha}, \ \forall t.
$$
It is obvious that $\Lambda_2^+\cap\Lambda_2^-\subset \{X_t> X_t^{\alpha},\forall t\}$. Hence
\begin{equation}\label{second}
\M\{X_t> X_t^{\alpha},\forall t\}\ge \M\{\Lambda_2^+\cap\Lambda_2^-\}=1-\alpha.
\end{equation}
Since the opposite of $\{X_t\le X_t^{\alpha},\forall t\}$ is $\{X_t\nleq X_t^{\alpha},\forall t\}$, on the basis of duality axiom, we can get
$$
\M\{X_t\le X_t^{\alpha},\forall t\}+\M\{X_t\nleq X_t^{\alpha},\forall t\}=1.
$$
Besides, $\{X_t> X_t^{\alpha},\forall t\}\subset\{X_t\nleq X_t^{\alpha},\forall t\}$ implies 
\begin{equation}\label{final}
\M\{X_t\le X_t^{\alpha},\forall t\}+\M\{X_t> X_t^{\alpha},\forall t\}\le 1.
\end{equation}

Thus, from the (\ref{first}),(\ref{second}) and (\ref{final}), it is evident that
\begin{equation*}
\begin{aligned}
&\M\{X_t\le X_t^{\alpha},\forall t\}=\alpha,\\
&\M\{X_t>X_t^{\alpha},\forall t\}=1-\alpha.
\end{aligned}
\end{equation*}
The theorem is proved.\\

Through the discussion of the above theorems, we can easily draw the following theorem.

\begin{Theorem}\label{Th:distribution1}
	Let $f\big(t,z_0,z_1,\cdots,z_{n-1}\big)$ be a function on $D=[0, a]\!\times\! \mathbb{D}_0\subset\mathbb{R}^{n\!+\!1}$ $(\mbox{where}\ \mathbb{D}_0\  \mbox{is an n-dimensional bounded region})$ satisfying local Lipschitz conditions in $z_0,z_1,\cdots,z_{n-1}$, and $g_1(t), g_2(t), \dots, g_m(t)$ be integrable functions with respect to $t$. Assume that $f\big(t,z_0,z_1,\cdots,z_{n-1}\big)$ is a monotonically increasing function with respect to $z_0,z_1,\cdots,$ $z_{n-2}$. 
	If $X_t$ and $X_t^{\alpha}$ are the solution of 
	\begin{equation}
	\frac{\mathrm{d}^n X_t}{\mathrm{d}t^n}=f\bigg(t,X_t,\cdots,\frac{\mathrm{d}^{n-1} X_t}{\mathrm{d}t^{n-1}}\bigg)+\sum_{i=1}^{m}g_i(t)\frac{\mathrm{d}C_{it}}{\mathrm{d}t}
	\end{equation}
	and 
    \begin{equation}
    \frac{\mathrm{d}^n X_t^{\alpha}}{\mathrm{d}t^n}=f\bigg(t,X_t^{\alpha},\cdots,\frac{\mathrm{d}^{n-1} X_t^{\alpha}}{\mathrm{d}t^{n-1}}\bigg)+\sum_{i=1}^{m}|g_i(t)|\Phi^{-1}(\alpha),
    \end{equation}
	respectively, where
	 $$
	\Phi^{-1}(\alpha)=\frac{\sqrt{3}}{\pi}\mathrm{ln}\frac{\alpha}{1-\alpha},\qquad 0<\alpha<1,
	$$ 
    then 
	 $X_t$ has an inverse uncertainty distribution
	\begin{equation}
	\Psi_t^{-1}(\alpha)=X_t^{\alpha}.
	\end{equation}
\end{Theorem}



\begin{Example}
	Suppose that $X_t$ and $X_t^{\alpha}$ are the solution of
	\begin{equation}\label{example3}
	\left\{
	\begin{aligned}
	&\frac{\mathrm{d}^2 X_t}{\mathrm{d}t^2}=2\frac{\mathrm{d} X_t}{\mathrm{d}t}+3X_t+e^{-t}\frac{\mathrm{d}C_{t}}{\mathrm{d}t}\\ &X_{(0)}=0\\
    &X'_{(0)}=0
	\end{aligned}
	\right.
	\end{equation}
    and
    \begin{equation}\label{ode}
    \left\{
    \begin{aligned}
    &\frac{\mathrm{d}^2 X^{\alpha}_t}{\mathrm{d}t^2}=2\frac{\mathrm{d} X^{\alpha}_t}{\mathrm{d}t}+3X^{\alpha}_t+e^{-t}\Phi^{-1}(\alpha)\\
    &X^{\alpha}_{(0)}=0\\
    &{X^{\alpha}_{(0)}}'=0
    \end{aligned}
    \right.
    \end{equation}
    where
    $$
    \Phi^{-1}(\alpha)=\frac{\sqrt{3}}{\pi}\mathrm{ln}\frac{\alpha}{1-\alpha},\qquad 0<\alpha<1,
    $$
    respectively.

    Note that $$2\frac{\mathrm{d} X^{\alpha}_t}{\mathrm{d}t}+3X^{\alpha}_t+e^{-t}\Phi^{-1}(\alpha)$$ is a monotonically increasing function with respect 
    to $X_t^{\alpha}$, so $X^{\alpha}_t$ is the $\alpha$-path of $X_t$.
    Solve this Eq.\eqref{ode} through Euler exponential function method and obtain the solution 
    $$
    X_t^{\alpha}=\frac{\sqrt{3}}{16\pi}\bigg(e^{3t}-e^{-t}-4te^{-t}\bigg)\cdot \ln\frac{\alpha}{1-\alpha},\quad 0<\alpha<1.
    $$ 
        
    \begin{figure}[htbp]
    	\centering
    	\includegraphics[scale=0.45]{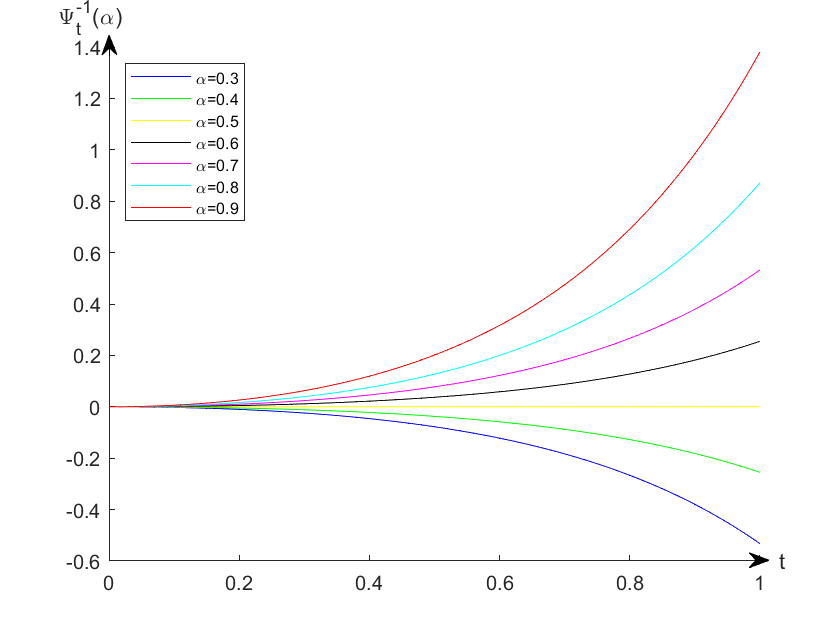}
    	\caption{Different values of $\alpha$ of $\Psi_t^{-1}(\alpha)$}
    	\label{graph3}
    \end{figure}

    According to Theorem \ref{Th:distribution1}, $X_t^{\alpha}$ is the $\alpha$-path of $X_t$ (shown in Fig.\ref{graph3}), and the inverse uncertainty distribution of $X_t$ is
    $$
    \Psi_t^{-1}(\alpha)=X_t^{\alpha}.
    $$

\end{Example}


\section{Residual}\label{residuals}
In general, the observations of an uncertain process are discrete points when the time intervals are not very short. In order to make a connection between HUDE and observations, we first define the residual of HUDE.
Let us consider a high-order uncertain differential equation 
\begin{equation*}
\frac{\mathrm{d}^n X_t}{\mathrm{d}t^n}=f\bigg(t,X_t,\frac{\mathrm{d} X_t}{\mathrm{d}t}, \cdots,\frac{\mathrm{d}^{n-1} X_t}{\mathrm{d}t^{n-1}}\bigg)+\sum_{i=1}^{m}g_i\bigg(t,X_t,\frac{\mathrm{d} X_t}{\mathrm{d}t},\cdots,\frac{\mathrm{d}^{n-1} X_t}{\mathrm{d}t^{n-1}}\bigg)\frac{\mathrm{d}C_{it}}{\mathrm{d}t}
\end{equation*}
where $f$ and $g_i\ (i=1,2,\cdots,m)$ are known continuous functions and $C_{1t},C_{2t},\cdots,$ $C_{mt}$ are independent Liu processes. Assume that 
\begin{equation}\label{data}
x_{t_j},\ x'_{t_j},\ ,\cdots,\ x^{(n\!-\!1)}_{t_j}
\end{equation}
are observations of $X_t,\frac{\mathrm{d} X_t}{\mathrm{d}t}, \cdots,\frac{\mathrm{d}^{n-1} X_t}{\mathrm{d}t^{n-1}}$ at $t_j\ (j=1,2,\cdots,l)$, where $t_1<t_2<\cdots<t_l$. Observation $x_{t_{l+1}}$ at time $t_{l+1}$ is obtained.
For convenience, all observed data is listed in Table \ref{table1}.

For any given index $j$ with $1\le j\le l$, we consider solving the updated high-order uncertain differential equation,
\begin{equation}\label{eq:updated}
\left\{
\begin{array}{l}
\vspace{1ex}
\displaystyle \frac{\mathrm{d}^n X_t}{\mathrm{d}t^n}\!\!=\!\!f\bigg(t,X_t,\frac{\mathrm{d} X_t}{\mathrm{d}t}, \cdots,\frac{\mathrm{d}^{n\!-\!1} X_t}{\mathrm{d}t^{n\!-\!1}}\bigg)\!+\!\sum\limits_{i=1}^{m}g_i\bigg(t,X_t,\frac{\mathrm{d} X_t}{\mathrm{d}t},\cdots,\frac{\mathrm{d}^{n\!-\!1} X_t}{\mathrm{d}t^{n\!-\!1}}\bigg)\frac{\mathrm{d}C_{it}}{\mathrm{d}t}\\
\vspace{1ex}
X_{t_{j}}=x_{t_{j}}\\
\displaystyle\frac{\mathrm{d} X_t}{\mathrm{d}t}\bigg|_{t=t_{j}}=x'_{t_{j}}\\
\qquad\quad\vdots\\
\displaystyle\frac{\mathrm{d}^{n-1} X_t}{\mathrm{d}t^{n-1}}\bigg|_{t=t_{j}}=x^{(n\!-\!1)}_{t_{j}}\\
\end{array}
\right.
\end{equation} 
where $x_{t_{j}},x'_{t_{j}},\cdots,x^{(n\!-\!1)}_{t_{j}}$ are observations at time $t_{j}$. 
The uncertainty distribution of $X_{t_{j+1}}$ initialized at $t_{j}$, denoted by $\Phi_{t_{j+1}}$, could be obtained by solving Eq.\eqref{eq:updated}.

\begin{table}[htp]
	\caption{Observed data at time $t_j$ ($1\le j \le l+1$)}\label{table1}
	\centering
	\begin{tabular}{ccccc}\hline
		$t_1$ & $t_2$ & $\cdots$ & $t_l$ & $t_{l+1}$\\ \hline
		$x_{t_1}$& $x_{t_2}$& $\cdots$& $x_{t_l}$ & $x_{t_{l+1}}$\\ 
		$x'_{t_1}$& $x'_{t_2}$& $\cdots$& $x'_{t_l}$\\ 
		$\vdots$& $\vdots$& $\vdots$& $\vdots$\\ 
		$x_{t_1}^{(n\!-\!1)}$& $x^{(n\!-\!1)}_{t_2}$& $\cdots$& $x^{(n\!-\!1)}_{t_l}$\\ \hline
	\end{tabular}
\end{table}

\begin{Proposition}\label{pro1}
  Let $\xi$ be a uncertain variable with regular uncertainty distribution $\Phi(x)$, then $\Phi(\xi)$ is a 
  linear uncertain variable $\mathcal{L}(0,1)$.
\end{Proposition}


For $1\le j \le l$, we know the uncertain distribution of $X_{t_{j+1}}$ and its observed value $x_{t_{j+1}}$.
Then the $j$th residual is defined as 
\begin{equation}\label{residual}
\varepsilon_j=\Phi_{t_{j+1}}(x_{t_{j+1}}).
\end{equation}
Then we have a total of $l$ residuals $$\varepsilon_1, ~\varepsilon_2, ~\cdots, ~\varepsilon_{l}.$$ 
By Proposition \ref{pro1},
all the residuals form a sample of linear uncertainty distribution $\mathcal{L}(0,1)$.


For most of high-order uncertain differential equations, it is not easy to obtain analytical solutions. Sometimes, it is even unpractical to calculate any analytical solutions. In this case, we need to use numerical methods to calculate the residual. The following algorithm is able to calculate the $j$th residual $\varepsilon_j$.\\
\ \\
\indent{\bf Algorithm 1: Numerical method for calculating residuals}\\
\indent{\bf Step 0:} Set $l=0$, $r=1$ and a precision $\delta=0.0001$.\\
\indent{\bf Step 1:} Set $\displaystyle\alpha= (l+r)/2$.\\
\indent{\bf Step 2:} Compute $X_{t_{j+1}}^{\alpha}$ of the updated Eq.\eqref{eq:updated} by Euler method.\\
\indent{\bf Step 3:} If $X_{t_{j+1}}^{\alpha}<x_{t_{j+1}}$, then $l= \alpha$. Otherwise, $r= \alpha$.\\
\indent{\bf Step 4:} If $|l-r|>\delta$, then go to Step 1.\\
\indent{\bf Step 5:} Output $\varepsilon_j= (l+r)/2$.\\
\par


\begin{Remark}
	However, in general, we can only obtain the observations of $X_t$ at different times, and it is difficult to get the value of $X'_t,\cdots,X^{(n-1)}_t$ at the corresponding time by observation. In this case, suppose 
	\begin{equation}\label{data3}
	x_{t_1},x_{t_2},\cdots,x_{t_{l+n-1}}
	\end{equation}
	are $l+n-1$ observations of $X_t$ at the times $t_1,t_2,\cdots,t_{l+n-1}$ with $t_1<t_2<\cdots<t_{l+n-1}$, respectively.
	Define
	\begin{equation}\label{es}
	x^{(v)}_{t_{j}}=\frac{x^{(v-1)}_{t_{j+1}}-x^{(v-1)}_{t_{j}}}{t_{j+1}-t_{j}}
	\end{equation}
	is the value of $X_{t}^{(v)}(v=1,2,\cdots,n-1)$ at the $t_j$ where $1\le j\le l+n-2$. Therefore, $l+n-1$ data are able to get the complete initial information of the first $l$ moments. All data is shown in Table \ref{table2}.
	
	\begin{table}[htp]
		\caption{Observed data at time $t_j$ ($1\le j \le l+n-1$)}\label{table2}
		\centering
		\begin{tabular}{cccccccc}\hline
			$t_1$ & $t_2$ & $\cdots$ & $t_l$ & $t_{l+1}$ & $\cdots$ & $t_{l+n}$ & $t_{l+n-1}$\\ \hline
			$x_{t_1}$& $x_{t_2}$& $\cdots$& $x_{t_l}$ & $x_{t_{l+1}}$  & $\cdots$ & $x_{t_{l+n}}$ & $x_{t_{l+n-1}}$\\ 
			$x'_{t_1}$& $x'_{t_2}$& $\cdots$& $x'_{t_l}$ &$x'_{t_{l+1}}$ & $\cdots$ & $x'_{t_{l+n}}$\\ 
			$\vdots$& $\vdots$&  & $\vdots$& $\vdots$ & &\\
			$x_{t_1}^{(n\!-\!2)}$& $x^{(n\!-\!2)}_{t_2}$& $\cdots$& $x^{(n\!-\!2)}_{t_l}$& $x^{(n\!-\!2)}_{t_{l+1}}$\\ 
			$x_{t_1}^{(n\!-\!1)}$& $x^{(n\!-\!1)}_{t_2}$& $\cdots$& $x^{(n\!-\!1)}_{t_l}$\\ \hline
		\end{tabular}
	\end{table}
Using difference \eqref{es} is the simplest form, in addition, derivative can be approximated by numerical differentiation or Lagrange interpolation process. For example, consider central difference
	$$
	x^{(v)}_{t_{j}}=\frac{x^{(v-1)}_{t_{j+1}}-x^{(v-1)}_{t_{j-1}}}{t_{j+1}-t_{j-1}}.
	$$
\end{Remark}

There are two main use of residuals: parameter estimation and hypothesis test, which will be discussed in the next two sections.



\section{Parameter estimation}\label{PE}

When we apply the model of HUDE to practical problems, the differential equations normally contain unknown 
parameters. In this section, we will 
estimate unknown parameters based on observed data and residuals.

Consider the following high-order uncertain differential equation
	\begin{equation}\label{eq:parameter hude}
	\frac{\mathrm{d}^n X_t}{\mathrm{d}t^n}=f\bigg(t,X_t, \cdots,\frac{\mathrm{d}^{n\!-\!1} X_t}{\mathrm{d}t^{n\!-\!1}};\bm{\theta}\bigg)+\sum_{i=1}^{m}g_i\bigg(t,X_t,\cdots,\frac{\mathrm{d}^{n\!-\!1} X_t}{\mathrm{d}t^{n\!-\!1}};\bm{\theta}\bigg)\frac{\mathrm{d}C_{it}}{\mathrm{d}t}
	\end{equation}
where $f$ and $g_i\ (i=1,2,\cdots,m)$ are known continuous functions, $C_{1t},C_{2t},\cdots,C_{mt}$ are independent Liu processes, and $\bm{\theta}(\bm{\theta}\in \mathcal{I}, \mathcal{I}\subseteq \mathbb{R}^p)$ is an unknown $p$-vector of parameters. Suppose the observed data is listed in Table \ref{table1}.

Then for every $j$ with $1\le j\le l$, the updated high-order uncertain differential equation
\begin{equation}\label{eq:updated-theta}
\left\{
\begin{array}{l}
\vspace{1ex}
\displaystyle\frac{\mathrm{d}^n X_t}{\mathrm{d}t^n}\!\!=\!\!f\bigg(t,X_t,\frac{\mathrm{d} X_t}{\mathrm{d}t}, \!\cdots\!,\frac{\mathrm{d}^{n\!-\!1} X_t}{\mathrm{d}t^{n\!-\!1}};\!\bm{\theta}\bigg)\!+\!\sum\limits_{i=1}^{m}g_i\bigg(t,X_t,\frac{\mathrm{d} X_t}{\mathrm{d}t},\!\cdots\!,\frac{\mathrm{d}^{n\!-\!1} X_t}{\mathrm{d}t^{n\!-\!1}};\!\bm{\theta}\bigg)\frac{\mathrm{d}C_{it}}{\mathrm{d}t}\\
\vspace{1ex}
X_{t_{j}}=x_{t_{j}}\vspace{1ex}\\
\displaystyle\frac{\mathrm{d} X_t}{\mathrm{d}t}\bigg|_{t=t_{j}}=x'_{t_{j}}\\
\qquad\quad\vdots\\
\displaystyle\frac{\mathrm{d}^{n-1} X_t}{\mathrm{d}t^{n-1}}\bigg|_{t=t_{j}}=x^{(n\!-\!1)}_{t_{j}}\\
\end{array}
\right.
\end{equation} 
contains unknown vector $\bm{\theta}$.
Then by Eq. \eqref{residual}, the $j$th residual is
\begin{equation}\label{residual1}
\varepsilon_j(\bm{\theta})=\Phi_{t_{j+1}}(x_{t_{j+1}}|\bm{\theta}).
\end{equation}
By Proposition \ref{pro1},
$$
\varepsilon_{1}(\bm{\theta}),\varepsilon_{2}(\bm{\theta}),\!\cdots\!,\varepsilon_{l}(\bm{\theta})
$$
form a sample of $\mathcal{L}(0,1)$. Based on the residuals, there are two methods to estimate the unknown $p$-vector $\bm{\theta}$. 

The first esitmation of $\bm{\theta}$ is the moment estimation. As the $k$th population moment of the linear uncertainty distribution $\mathcal{L}(0,1)$ is
$$
\frac{1}{k+1}.
$$
According to the principle that the $k$th sample moment is equal to the $k$th population moment, the moment estimate $\bm{\theta}$ should resolve the system of equations
\begin{equation}\label{ude estimation}
\frac{1}{l}\sum_{j=1}^{l}\varepsilon^{k}_{j}(\bm{\theta})=\frac{1}{k+1}, \quad k=1,2,\cdots,p,
\end{equation}
where $p$ is the dimension of $\bm{\theta}$.
Thus, the solution of the 
following minimization problem,

\begin{equation}\label{problem}
\left\{
\begin{array}{ll}
\vspace{1mm}
\min\limits_{\bm{\theta}}\sum\limits_{k=1}^{p}\bigg(\displaystyle\frac{1}{l}\sum\limits_{j=1}^{l}\varepsilon^{k}_{j}(\bm{\theta})-\frac{1}{k+1}\bigg) ^2\\
\rm{subject\ to:}\\
\ \qquad\qquad \bm{\theta}\in \mathcal{I}
\end{array}
\right.
\end{equation}
is the moment estimation of $\bm{\theta}$. This minimization problem \eqref{problem} could be solved by MATLAB\footnote{MATLAB R2021a, 9.10.0.1602886, maci64, Optimization Toolbox, ``fminsearch" function.}.


\begin{Remark}
	The optimal value of the objective function in Eq \eqref{problem} should be very close to zero. In the actual calculation using MATLAB, it is generally considered that the value of objective function should be less than $10^{-10}$. Otherwise, we can assume that data do not fit the proposed uncertain differential equation. 
\end{Remark}

The other estimation is the maximum likelihood estimation proposed by \cite{LiuandLiu2023-1}. 
Given a detection level $\alpha$, the maximum likelihood estimation of $\bm{\theta}$ is the solution of the following system of equations,
\begin{equation*}
	\begin{cases}
	\displaystyle\varepsilon'_{i^{*}(\bm{\theta})}(\bm{\theta})
	=\frac{\alpha}{2}& \\
	\displaystyle\varepsilon'_{i^{*}(\bm{\theta})+\lceil l(1-\alpha) \rceil-1}(\bm{\theta})
	=1-\frac{\alpha}{2}& \\
	i^{*}(\bm{\theta})=\mathop{\arg\min}\limits_{1\le i\le l-\lceil l(1-\alpha) \rceil+2}\varepsilon'_{i+\lceil l(1-\alpha) \rceil-1}(\bm{\theta})-\varepsilon'_{i}(\bm{\theta}),
\end{cases}
\end{equation*}
where 
\[
\{\varepsilon'_{1}(\bm{\theta}),\varepsilon'_{2}(\bm{\theta}),\cdots,\varepsilon'_{l}(\bm{\theta})\}
\] 
with
\[
\varepsilon'_{1}(\bm{\theta})\le\varepsilon'_{2}(\bm{\theta}),\cdots\le\varepsilon'_{l}(\bm{\theta})
\]
is a rearrangement of 
\[
\{\varepsilon_1(\bm{\theta}),\varepsilon_2(\bm{\theta}),\cdots,\varepsilon_{l}(\bm{\theta})\}.
\]



\section{Uncertain hypothesis test}\label{UHT}
After the unknown parameters have been estimated using the method in Sect.\ref{PE}, it is crucial to test whether a HUDE is a good fit to the observed datas. 
Here we will use the uncertain hypothesis test introdueced by \cite{LiuandYe2023} to evaluate the fitness of HUDE. 

Consider an high-order uncertain differential equation
\begin{equation}\label{eq:1}
\frac{\mathrm{d}^n X_t}{\mathrm{d}t^n}=f\bigg(t,X_t, \cdots,\frac{\mathrm{d}^{n\!-\!1} X_t}{\mathrm{d}t^{n\!-\!1}};\bm{\theta}\bigg)+\sum_{i=1}^{m}g_i\bigg(t,X_t,\cdots,\frac{\mathrm{d}^{n\!-\!1} X_t}{\mathrm{d}t^{n\!-\!1}};\bm{\theta}\bigg)\frac{\mathrm{d}C_{it}}{\mathrm{d}t}
\end{equation}
where $f$ and $g_i\ (i=1,2,\cdots,m)$ are known continuous functions and $C_{1t},C_{2t},\cdots,$ $C_{mt}$ are independent Liu processes but $\bm{\theta}(\bm{\theta}\in \mathcal{I}, \mathcal{I}\subseteq \mathbb{R}^p)$ is an unknown $p$-vector of parameters. Assume that 
\begin{equation}\label{data2}
x_{t_1},x_{t_2},\cdots,x_{t_l}
\end{equation}
are $l$ observations of $X_t$ at the times $t_1,t_2,\cdots,t_l$ with $t_1<t_2<\cdots<t_l$, where $l>n$, respectively.
For each index $j(1\le j\le l-n+1)$, we solve the updated high-order uncertain differential equation  
\begin{equation}\label{eq:uht}
\left\{
\begin{array}{l}
\vspace{1ex}
\displaystyle\frac{\mathrm{d}^n X_t}{\mathrm{d}t^n}\!=\!f\bigg(t,X_t,\frac{\mathrm{d} X_t}{\mathrm{d}t},\! \cdots\!,\frac{\mathrm{d}^{n-1} X_t}{\mathrm{d}t^{n-1}};\bm{\theta}\bigg)\!+\!\sum\limits_{i=1}^{m}g_i\bigg(t,X_t,\frac{\mathrm{d} X_t}{\mathrm{d}t},\!\cdots\!,\frac{\mathrm{d}^{n\!-\!1} X_t}{\mathrm{d}t^{n\!-\!1}};\bm{\theta}\bigg)\frac{\mathrm{d}C_{it}}{\mathrm{d}t}\\
\vspace{1ex}
X_{t_{j}}=x_{t_{j}}\\
\displaystyle\frac{\mathrm{d} X_t}{\mathrm{d}t}\bigg|_{t=t_{j}}=x'_{t_{j}}\\
\qquad\quad\vdots\\
\displaystyle\frac{\mathrm{d}^{n-1} X_t}{\mathrm{d}t^{n-1}}\bigg|_{t=t_{j}}=x^{(n-1)}_{t_{j}}\\
\end{array}
\right.
\end{equation}
where $x_{t_{j}},x'_{t_{j}},\cdots,x^{(n\!-\!1)}_{t_{j}}$ are $n$ new initial values according to the form of difference \eqref{es} at the new initial time $t_{j}$ with $1\le j\le l-n+1$, respectively.

For any given $\bm{\theta}$, on the basis of Sec. \ref{PE}, we can produce $l-n+1$ residuals 
$$
\varepsilon_{1}(\bm{\theta}),\varepsilon_{2}(\bm{\theta}),\!\cdots\!,\varepsilon_{l\!-\!n+1}(\bm{\theta})
$$
of Eq.\eqref{eq:uht} corresponding to the observed data \eqref{data2}.

If the high-order uncertain differential Eq.\eqref{eq:1} does fit the observed data \eqref{data2}, then
$$
\varepsilon_{1},\varepsilon_{2},\!\cdots\!,\varepsilon_{l\!-\!n+1}\sim \mathcal{L}(0,1).
$$
That is, testing whether a high-order uncertain differential equation fits the observed data well is equivalent to testing whether these residuals $\varepsilon_{1},\varepsilon_{2},\!\cdots\!,\varepsilon_{l\!-\!n+1}$ fit the uncertainty distribution $\mathcal{L}(0,1)$, i.e.
$$
\varepsilon_{1},\varepsilon_{2},\!\cdots\!,\varepsilon_{l\!-\!n+1}\sim \mathcal{L}(0,1).
$$

Given a significance level $\alpha\ (\text{e.g.}\  0.05)$, the test is 
\begin{equation*}
\begin{aligned}
W=\bigg\{(\varepsilon_{1}, \varepsilon_{2}, \cdots, \varepsilon_{l\!-\!n+1}):&{\text{there are at least $\alpha$ of indexes j's with }} 1\le j\le l\!-\!n\!+\!1 \\
&{\text{ such that }}\varepsilon_j<\frac{\alpha}{2} {\text{ or }} \varepsilon_j>1-\frac{\alpha}{2}\bigg\}.
\end{aligned}
\end{equation*}
If the vector of the $l-n+1$ residuals $\varepsilon_{1}, \varepsilon_{2}, \cdots, \varepsilon_{l\!-\!n+1}$ belongs to the test $W$, i.e.,
$$
(\varepsilon_{1}, \varepsilon_{2}, \cdots, \varepsilon_{l\!-\!n+1})\in W,
$$
then Eq.\eqref{eq:1} is not a good fit to the observed data \eqref{data2}.

\noindent If
$$
(\varepsilon_{1}, \varepsilon_{2}, \cdots, \varepsilon_{l\!-\!n+1})\notin W,
$$
then Eq.\eqref{eq:1} is a good fit to the observed data \eqref{data2}.


\section{An application in nuclear reactor kinetics}\label{nuclear reactor}
Nuclear safety has emerged as a crucial consideration for many nations in their pursuit of peaceful development. Despite the universal commitment to nuclear disarmament and non-proliferation, research in the nuclear industry continues unabated.

The nuclear reactor kinetics equations have been studied and modeled by many scholars(\cite{AllenandHayes2005}; \cite{Allen1999,Allen2004}). The nuclear reactor kinetics is also known as neutron population kinetics which studies the dynamic change of neutron population in reactor. Neutron population determines the change of power level over time and are affected by the control rod position and other factors(\cite{AllenandHayes2005}). 


Due to the specific operational complexities of nuclear reactors, research into refining nuclear reactor models is imperative. Furthermore, in cases involving newly discovered radioactive isotopes or neutrons lacking comprehensive data, experts often rely on past experience to give credibility to analyze their properties. The current nuclear reactor models do not adequately account for such uncertainties.

\subsection{Nuclear reactor kinetics driven by Liu process}\label{deduce}
\noindent In this part, the quantitative relationship between neutron population and time variation in a nuclear reactor is modeled. In a real nuclear reactor, if the delayed neutrons are not taken into account, the neutron population increases so rapidly in the supercritical state that the reactor cannot be controlled, which is extremely dangerous for a nuclear reactor. A delayed neutron was produced by a delayed neutron precursor. Every delayed neutron precursor likes Br-87, Uranium-235 is produced at the instant of fission and releases delayed neutrons after a slight delay. In order to better establish the nuclear reactor kinetics model, the nuclear reactor is supposed to be large enough that energy and space effects can be ignored.

\begin{figure}[htbp]
	\centering
	\includegraphics[scale=0.45]{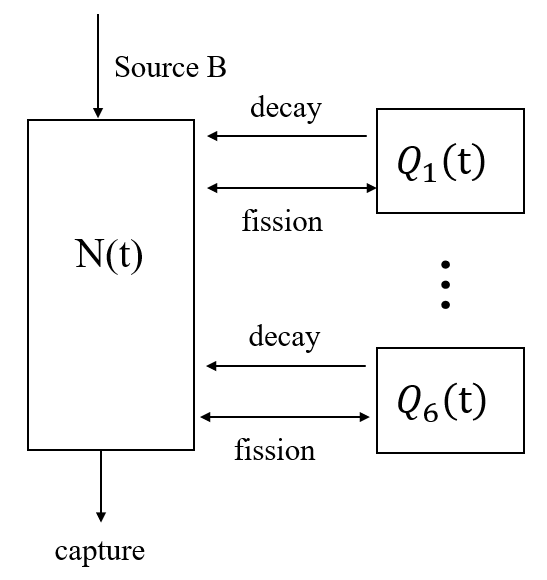}
	\caption{Dynamical system of a nuclear reactor}
	\label{graph1}
\end{figure}
Normally, there is more than one delayed neutron precursors in a reactor, and nuclear physicists divide these delayed neutron precursors into roughly six groups based on average lifetime and decay constant. Therefore, it is assumed that there are six groups of delayed neutron precursors in this reactor. Fig.\ref{graph1} depicts the dynamic processes of the neutron concentrations and precursors denisty.

In this figure, an extraneous neutron source $B$ transports the neutron at a constant rate to the nuclear reactor. $N(t)$ is the neutron population with respect to time $t$ and the variable $Q_{i}(t)$ is $i$-th delayed neutron precursor density over time. The following parameters are important to describe the dynamical system.
Assume that $k$ is effective multiplication constant of neutron which is calculated as neutrons population in the new generation divided by neutrons population in the old generation
and neutron lifetime is $l$. Consider that $\beta_i$ is the fraction of $i$-th delayed neutron precursor and $\beta=\sum\limits_{i=1}\limits^{6}\beta_i$ is the delayed neutron fraction. And $\lambda_i$ is the decay constant of fission product isotope $Q_{i}$.  

The entire dynamic process of nuclear reactor is briefly described as follows. It is not difficult to draw that if the initial value of neutrons is $N_t$, neutron population through a nuclear fission is $kN_t$ where delayed neutron density is $k\beta N_t$ and the number of new neutrons is $k(1-\beta)N_t$. 
At first, discuss the change of $i$-th delayed neutron precursors density $Q_i(t)$ with respect to time interval $\Delta t$. When a fission occurs, this reactor produces delayed neutron precursors $Q_i$ with a certain rate $\frac{k\beta_iN_t}{l}$ . Another situation, $i$-th delayed neutron precursor decays into a neutron with decay constant $\lambda_i$. Then the change of $Q_i$ is

\begin{equation}\label{eq:Q}
\Delta Q_{i}=\bigg(-\lambda_i Q_i+\frac{k\beta_iN_t}{l}\bigg)\Delta t.
\end{equation}

Next, consider the change of the neutron population $N(t)$ with respect to time interval $\Delta t$. There are three situations about the change of the neutron population which are source, transformation and born. The source event represents that an extraneous neutron source delivers $B$ neutrons to this reactor at a constant rate, so derive $B\Delta t$ neutrons at time interval $\Delta t$. The transformation event represents that $i$-th delayed neutron precursor decays into a neutron with decay constant $\lambda_i$. The born event represents neutrons produced when a fission occurs with rate $\frac{k(1-\beta)-1}{l}N_t$.
Given the above analysis, we can get the change of $N(t)$:
\begin{equation}\label{eq:N(t)}
\Delta N_t=\bigg(B+\sum_{i=1}^{6}\lambda_i Q_i+\frac{k(1-\beta)-1}{l} N_t\bigg)\Delta t.
\end{equation}

According to Eq.\eqref{eq:N(t)} and Eq.\eqref{eq:Q}, we can derive the following formula about $N(t)$ and $Q_i(t)$

	\begin{equation}\label{eq:kinetics6}
	\left\{
	\begin{array}{l}
	\displaystyle\frac{\mathrm{d} N_t}{\mathrm{d}t}=B+\sum\limits_{i=1}^{6}\lambda_i Q_i+\frac{k(1-\beta)-1}{l} N_t\\
	
	\displaystyle\frac{\mathrm{d} Q_i}{\mathrm{d}t}=-\lambda_i Q_i+\frac{k\beta_i}{l}N_t,\quad i=1,2,\cdots,6.
	\end{array}
	\right.
	\end{equation}

The six groups of delayed neutron dynamics Eqs. \eqref{eq:kinetics6} are still complicated in form.
Considering a relatively stable situation, assume just one group of delayed neutron precursor in this reactor. We can get $\lambda=\lambda_1$ and $\beta=\beta_1$. In this way, the seven equations are reduced to two equations, so the equations are a little easier to solve.

In the beginning, the neutron population is very low in the reactor. In order to accelerate the start-up speed, it is necessary to add an extraneous neutron source B. As the reaction progresses, neutron population is much greater than the source neutron concentration, so we can regard as $B=0$. Thus we can get the equation of just one delayed neutron precursor:

	\begin{equation}\label{eq:kinetics1}
	\left\{
	\begin{array}{l}
	\vspace{5pt}
	\displaystyle\frac{\mathrm{d} N_t}{\mathrm{d}t}=\lambda Q_t+\frac{k(1-\beta)-1}{l} N_t\\
	
	\displaystyle\frac{\mathrm{d} Q_t}{\mathrm{d}t}=-\lambda Q_t+\frac{k\beta}{l}N_t.
	\end{array}
	\right.
	\end{equation}

The nuclear reactor kinetics of one delayed neutron can roughly reflect the fast transient of neutron population change and the slow transient with stable period.
Since we care more about the neutron population and the properties of this system of Eq.(\ref{eq:kinetics1}), we first consider transforming it into a second-order equation
\begin{equation}\label{eq:ode-2}
\frac{\mathrm{d}^2 N_t}{\mathrm{d}t^2}=\bigg(\frac{k(1-\beta)-1}{l}-\lambda\bigg)\frac{\mathrm{d}N_t}{\mathrm{d}t}+\lambda\cdot\frac{k-1}{l}N_t.
\end{equation}

The actual expression of the above Eq.(\ref{eq:ode-2}) is that the neutron population varies with time in the same way. However, such a situation will only occur in ideal conditions, or in computer simulations. 
Similar to the situation of start-up, shutdown in reactor or discovery of unknown new radioactive elements, we must fully consider the impact of unknown situations on the reactor to ensure the safe operation of the reactor.
Since in the actual reaction, $\beta$ and $\lambda$ will be disturbed by environmental factors and will not be a stable constant, that is
$$
\beta(t)=\beta+\sigma_1\cdot ``\text{noise}"
$$
$$
\lambda(t)=\lambda+\sigma_2\cdot ``\text{noise}".
$$
Based on uncertainty theory, the ``noise" can be consider as a normal uncertain variable $\mathcal{N}(0,1)$. In other words,
$$
\frac{\mathrm dC_{1t}}{\mathrm dt}=\frac{C_{1,t+\Delta t}-C_{1t}}{\Delta t}
$$
and
$$
\frac{\mathrm dC_{2t}}{\mathrm dt}=\frac{C_{2,t+\Delta t}-C_{2t}}{\Delta t}
$$
where $C_{1t}$ and $C_{2t}$ are Liu processes, so we get
$$
\left\{
\begin{array}{l}
\displaystyle\beta(t)=\beta+\sigma_1\cdot \frac{\mathrm dC_{1t}}{\mathrm dt}\vspace{5pt}\\
\displaystyle\lambda(t)=\lambda+\sigma_2\cdot \frac{\mathrm dC_{2t}}{\mathrm dt}
\end{array}
\right.
$$
where $\sigma_1$ and $\sigma_2$ are nonnegative constants, and they represent the noise levels.
Consequently, we derive the second-order uncertain differential equation of nuclear reactor kinetics, i.e.
	\begin{equation}\label{eq:UDE}
	\left\{
	\begin{array}{l}
	\vspace{1ex}	
	\displaystyle\frac{\mathrm{d}^2 N_t}{\mathrm{d}t^2}\!=\!\Big(\frac{k(1\!-\!\beta)\!-\!1}{l}\!-\!\lambda\Big)\frac{\mathrm{d}N_t}{\mathrm{d}t}\!+\!\lambda\frac{k\!-\!\!1}{l}N_t\!-\!\frac{k}{l}\sigma_1\frac{\mathrm{d}N_t}{\mathrm{d}t}\frac{\mathrm{d}C_{1t}}{\mathrm{d}t}\\
	\ \ \ \ \ \ \ \ \ \ \ \ \  \displaystyle+\sigma_2\big(\frac{k\!-\!1}{l}  N_t\!-\!\frac{\mathrm{d}N_t}{\mathrm{d}t}\big)\frac{\mathrm{d}C_{2t}}{\mathrm{d}t}\\
	\vspace{1ex}	
	N(t_0)=n_0\\	
	\displaystyle\frac{\mathrm{d}N_t}{\mathrm{d}t}\bigg|_{t=t_0}=n_0'.
	\end{array}
	\right.
	\end{equation}

\subsection{Parameter estimation and hypothesis test}\label{Parameter Estimation}

Although we deduced the uncertain nuclear reactor kinetics Eq.\eqref{eq:UDE}, only in connection with the actual data can we reflect whether this equation is reasonable or not.
In this part, we will utilize experimental data to estimate the unknown parameters in this equation. To safeguard real data, the utilized data are based on specific parameters of an actual reactor.

Regard $\lambda=0.0785/sec, \beta=0.0065, k=1.001, l=10^{-4}sec $ as a numerical example, which are the data involved in the thermal fission in a nuclear reactor about uranium-$235$ fuel.

Consider the updated uncertain differential equation of nuclear reactor kinetics
\begin{equation}\label{eq:data}
	\left\{
	\begin{array}{l}
		\vspace{1ex}
		\displaystyle\frac{\mathrm{d}^2 N_t}{\mathrm{d}t^2}\!=\!-55.1435\frac{\mathrm{d}N_t}{\mathrm{d}t}\!+\!0.785N_t
		\!-\!10010 \sigma_1\frac{\mathrm{d}N_t}{\mathrm{d}t}\frac{\mathrm{d}C_{1t}}{\mathrm{d}t}\\
		\ \ \ \ \ \ \ \ \ \ \ \ \displaystyle+\sigma_2\bigg(10 N_t\!-\!\frac{\mathrm{d}N_t}{\mathrm{d}t}\bigg)\frac{\mathrm{d}C_{2t}}{\mathrm{d}t}\\
		\vspace{1ex}
		N_{t_{j}}=x_{t_{j}}\vspace{1ex}\\
		\displaystyle\frac{\mathrm{d} N_t}{\mathrm{d}t}\bigg|_{t=t_{j}}=x'_{t_{j}}\\
	\end{array}
	\right.
\end{equation} 
where $x_{t_j}$, $x'_{t_j}$ are observations at time $t_j$, and $\sigma_1$ and $\sigma_2$ are unknown parameters to be estimated.

\begin{figure}[htbp]
	\centering
	\includegraphics[scale=0.2]{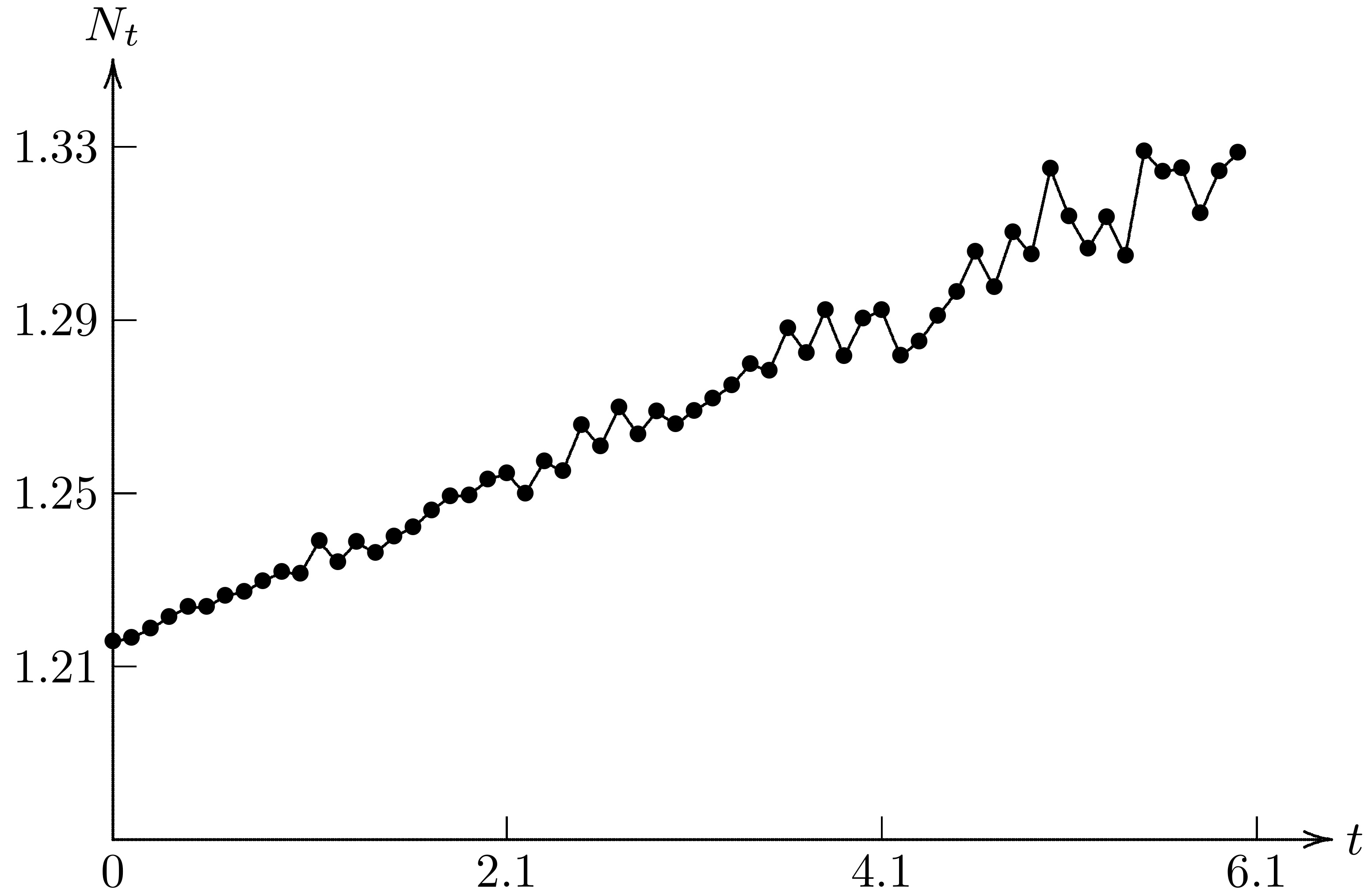}
	\caption{Neutron population from $t=0$ to $t=6$ with initial value $n_0=1.2157$}
	\label{graph2}
\end{figure}

Table \ref{61observation} shows $61$ observed data of this unclear reactor based on the values given above. And the fluctuation of the  neutron population over time is shown in the Fig.\ref{graph2}.

Next, we will apply residuals and moment estimation to estimate the unknown parameters. According to the Algorithm 1, we can produce 60 residuals 
$$
\varepsilon_1(\sigma_1,\sigma_2),\ \varepsilon_2(\sigma_1,\sigma_2),\ \cdots,\ 
\varepsilon_{60}(\sigma_1,\sigma_2)
$$
for any given parameters $\sigma_1$ and $\sigma_2$.

\begin{table}
	\caption{61 experimental data of nuclear reactor}\label{61observation}
	\centering
	\resizebox{\textwidth}{25mm}{
		\begin{tabular}{cccccccccccc}
			\hline
			$t$&$N_t$&$t$&$N_t$&$t$&$N_t$&$t$&$N_t$&$t$&$N_t$&$t$&$N_t$\\
			\hline
			$0$& $1.2157$& $1$& $1.2313$&	$2$&	$1.2530$&	$3$&	$1.2658$&	$4$&	$1.2902$& $5$&	$1.3249$\\
			$0.1$& $1.2165$& $1.1$& $1.2388$&	$2.1$&	$1.2545$&	$3.1$&	$1.2689$&	$4.1$&	$1.2922$& $5.1$&	$1.3138$\\
			$0.2$& $1.2186$& $1.2$& $1.2339$&	$2.2$&	$1.2498$&	$3.2$&	$1.2717$&	$4.2$&	$1.2817$& $5.2$&	$1.3064$\\
			$0.3$& $1.2213$& $1.3$& $1.2386$&	$2.3$&	$1.2572$&	$3.3$&	$1.2748$&	$4.3$&	$1.2849$& $5.3$&	$1.3136$\\
			$0.4$& $1.2236$& $1.4$& $1.2361$&	$2.4$&	$1.2550$&	$3.4$&	$1.2797$&	$4.4$&	$1.2908$&   $5.4$&	$1.3047$\\
			$0.5$& $1.2236$& $1.5$& $1.2399$&	$2.5$&	$1.2656$&	$3.5$&	$1.2782$&	$4.5$&	$1.2964$& $5.5$&	$1.3288$\\
			$0.6$& $1.2262$& $1.6$& $1.2420$&	$2.6$&	$1.2607$&	$3.6$&	$1.2880$&	$4.6$&	$1.3057$&   $5.6$&	$1.3241$\\
			$0.7$& $1.2271$& $1.7$& $1.2459$&	$2.7$&	$1.2697$&	$3.7$&	$1.2823$&	$4.7$&	$1.2975$&   $5.7$&	$1.3250$\\
			$0.8$& $1.2296$& $1.8$& $1.2492$&	$2.8$&	$1.2635$&	$3.8$&	$1.2922$&	$4.8$&	$1.3101$&   $5.8$&	$1.3145$\\
			$0.9$& $1.2317$& $1.9$& $1.2494$&	$2.9$&	$1.2688$&	$3.9$&	$1.2815$&	$4.9$&	$1.3050$&   $5.9$&	$1.3243$\\
			&&&&&&&&&&$6$&$1.3285$\\
			\hline
	\end{tabular}}
\end{table}
According to Eq.\eqref{problem},   
 the generalized moment estimation ($\hat{\sigma_1},\hat{\sigma_2}$) for ($\sigma_1,\sigma_2$) is the optimal solution of 
 \begin{equation}\label{estimation}
 	\left\{
 	\begin{array}{ll}
 		\vspace{1mm}
        \displaystyle\min\limits_{\sigma_1,\sigma_2}\sum\limits_{q=1}^{2}\bigg(\frac{1}{60}\sum_{i=1}^{60}\varepsilon_{i}^{q}(\sigma_1,\sigma_2)-\frac{1}{q+1}\bigg)^2\\
 		\rm{subject\ to:}\\
 		\ \qquad\qquad \sigma_1,\sigma_2\in (0,1).
 	\end{array}
 	\right.
 \end{equation}
Solving the above minimization problem \eqref{estimation} by MATLAB, we can obtain
$$
\hat{\sigma_1}=0.000143, \qquad \hat{\sigma_2}=0.296798,
$$
and the minimum value of Eq.(\ref{estimation}) is 
$$
\min_{\sigma_1,\sigma_2}\sum_{q=1}^{2}\bigg(\frac{1}{60}\sum_{i=1}^{60}\varepsilon_{i}^{q}(\sigma_1,\sigma_2)-\frac{1}{q+1}\bigg)^2 = 8.25350\times 10^{-12}.
$$

Thus we obtain uncertain differential equation of this nuclear reactor 
\begin{equation}\label{num}
	\left\{
	\begin{array}{ll}
		\vspace{2ex}	
		\displaystyle\frac{\mathrm{d}^2 N_t}{\mathrm{d}t^2}=\displaystyle-55.1435\frac{\mathrm{d}N_t}{\mathrm{d}t}\!+\!0.785N_t
		\!-\!1.43143 \frac{\mathrm{d}N_t}{\mathrm{d}t}\cdot\frac{\mathrm{d}C_{1t}}{\mathrm{d}t}\!\\
		\ \ \ \ \ \ \ \ \ \ \ \ \ \displaystyle+\bigg(2.96798N_t\!-\!0.296798\frac{\mathrm{d}N_t}{\mathrm{d}t}\bigg)\frac{\mathrm{d}C_{2t}}{\mathrm{d}t}\\
		\vspace{2ex}	
		N(0)=1.2157\\	
		N'(0)=0.008.
	\end{array}
	\right.
\end{equation}
Moreover, we can also get 60 residuals $\varepsilon_{1},\cdots,\varepsilon_{60}$ of uncertain nuclear reactor kinetics Eq.\eqref{num} as shown in Table~\ref{60residuals}. On the basis of hypothesis test, we need to test whether these 60 residuals $\varepsilon_{1},\cdots,\varepsilon_{60}$ have a good fit to the linear uncertainty distribution $\mathcal{L}(0,1)$. The test method is as follows. 

\begin{table}
	\caption{60 residuals}\label{60residuals}
	\centering
	\resizebox{\textwidth}{23mm}{
		\begin{tabular}{cccccccccccc}
			\hline
			$j$&$\varepsilon_j$&$j$&$\varepsilon_j$&$j$&$\varepsilon_j$&$j$&$\varepsilon_j$&$j$&$\varepsilon_j$&$j$&$\varepsilon_j$\\
			\hline
			$1$& $0.4373$& $11$& $0.8227$&	$21$&	$0.4787$&	$31$&	$0.5851$&	$41$&	$0.5070$& $51$&	$0.0317$\\
			$2$& $0.5238$& $12$& $0.1370$&	$22$&	$0.1454$&	$32$&	$0.5664$&	$42$&	$0.0342$& $52$&	$0.0809$\\
			$3$& $0.5632$& $13$& $0.6916$&	$23$&	$0.8153$&	$33$&	$0.5845$&	$43$&	$0.5906$& $53$&	$0.7937$\\
			$4$& $0.5436$& $14$& $0.2351$&	$24$&	$0.2537$&	$34$&	$0.6861$&	$44$&	$0.7404$& $54$&	$0.0540$\\
			$5$& $0.3800$& $15$& $0.6358$&	$25$&	$0.9065$&	$35$&	$0.2970$&	$45$&	$0.7223$&   $55$&	$0.9952$\\
			$6$& $0.5616$& $16$& $0.5177$&	$26$&	$0.1407$&	$36$&	$0.8839$&	$46$&	$0.8696$& $56$&	$0.1549$\\
			$7$& $0.4387$& $17$& $0.6417$&	$27$&	$0.8673$&	$37$&	$0.1187$&	$47$&	$0.0638$&   $57$&	$0.4357$\\
			$8$& $0.5498$& $18$& $0.6019$&	$28$&	$0.1011$&	$38$&	$0.8856$&	$48$&	$0.9360$&   $58$&	$0.0368$\\
			$9$& $0.5226$& $19$& $0.3940$&	$29$&	$0.7186$&	$39$&	$0.0329$&	$49$&	$0.1401$&   $59$&	$0.8784$\\
			$10$& $0.3566$& $20$& $0.6192$&	$30$&	$0.2134$&	$40$&	$0.8542$&	$50$&	$0.9875$&   $60$&	$0.6384$\\
			\hline
	\end{tabular}}
\end{table}

\begin{figure}[htbp]
	\centering
	\includegraphics[scale=0.2]{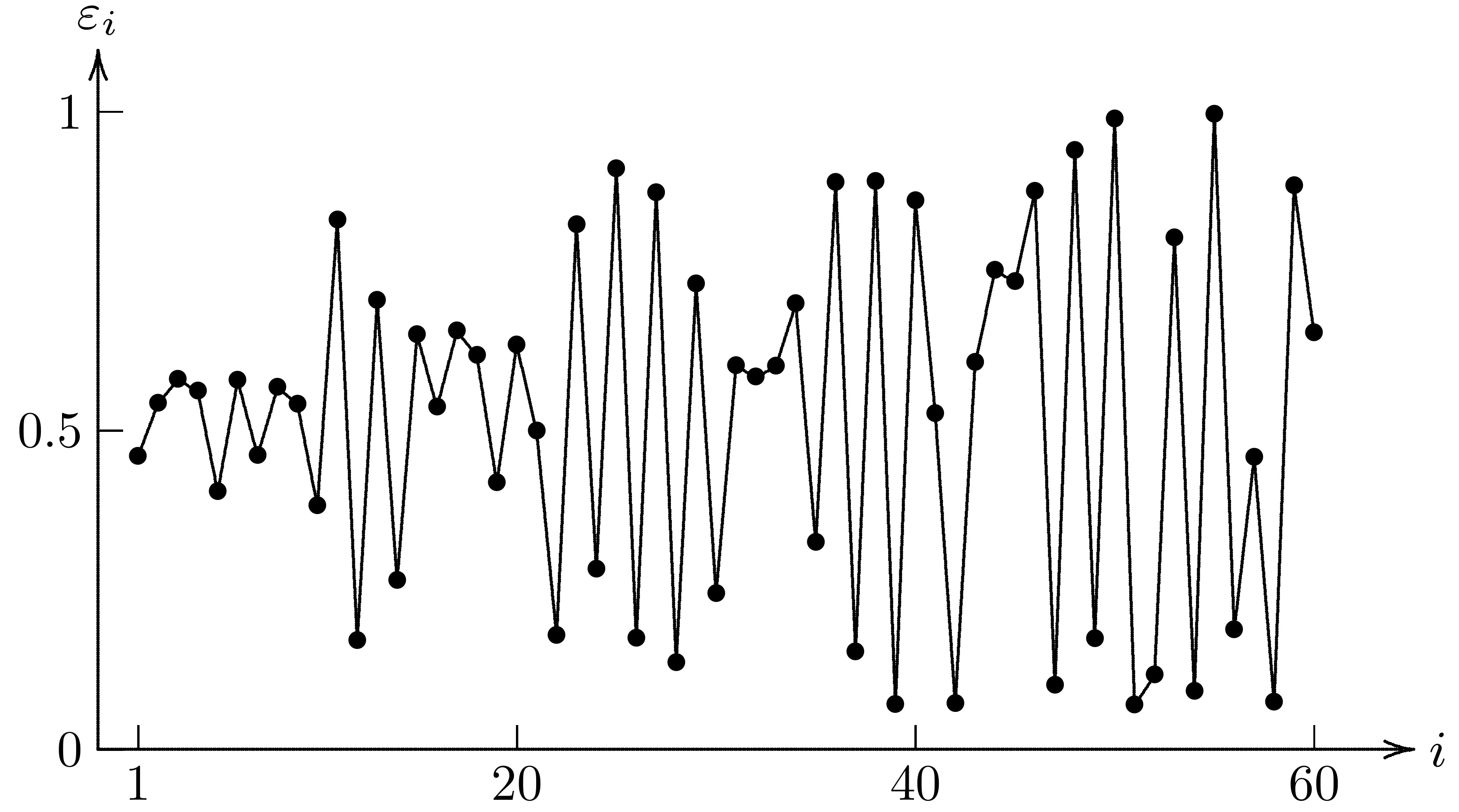}
	\caption{Residual plot of hypothesis test}
	\label{graph4}
\end{figure}

Given a significance level $\alpha=0.05$, and due to $\alpha\times 60=3$, the test is 
\begin{equation*}
	\begin{aligned}
		W=\{(\varepsilon_{1},\cdots,\varepsilon_{60}):&{\text{there are at least 3 of index j's with }} 1\le j\le 60 \\
		&{\text{ such that }}\varepsilon_j<0.025 {\text{ or }} \varepsilon_j>0.975\}.
	\end{aligned}
\end{equation*}
From Table \ref{60residuals} and Fig.\ref{graph4}, it is clearly that $\varepsilon_{50}$, $\varepsilon_{55}$ are the only two residuals not in $[0.025, 0.975]$. Therefore, $(\varepsilon_{1},\cdots,\varepsilon_{60})\notin W$. In other words, Eq.(\ref{num}) really have a good fit with the observed data.

Why do we use uncertain differential equations to describe nuclear reactor kinetics? The reason is as follows.
When we divide the 60 residuals into 
$$
(\varepsilon_{1},\varepsilon_{2},\cdots,\varepsilon_{21})\ \mbox{and}\  (\varepsilon_{44},\varepsilon_{45},\cdots,\varepsilon_{60}),
$$
the two-sample Kolmogorov-Smirnov test showed that the above two parts from the residuals do not come from the same population via the function "kstest2" in Matlab. Thus the residuals $\varepsilon_{1},\cdots,\varepsilon_{60}$ are not white noise in the sense of probability theory.

Next, consider the $\alpha$-path of the uncertain nuclear reactor kinetics Eq.\eqref{num}. $N_t^{\alpha}$ is the solution of 
	\begin{equation}\label{num1}
	\left\{
	\begin{array}{l}
		\displaystyle\frac{\mathrm{d}^2 N_t^{\alpha}}{\mathrm{d}t^2}=\displaystyle-55.1435\frac{\mathrm{d}N_t^{\alpha}}{\mathrm{d}t}\!+\!0.785N_t^{\alpha}
		\!+\bigg|-\!1.43143 \frac{\mathrm{d}N_t^{\alpha}}{\mathrm{d}t}\bigg|\Phi^{-1}(\alpha)\\
		\ \ \ \ \ \ \ \ \ \ \ \ \ \  \displaystyle+\bigg|2.96798N_t^{\alpha}\!-\!0.296798\frac{\mathrm{d}N_t^{\alpha}}{\mathrm{d}t}\bigg|\Phi^{-1}(\alpha)\\
		\vspace{2ex}	
		N(0)=1.2157\\	
		N'(0)=0.008.
	\end{array}
	\right.
    \end{equation}

In an actual controlled reactor, as the nuclear reaction goes on, the value of $\frac{\mathrm{d}N_t}{\mathrm{d}t}/N_t$ is very small, so that $N_t$ is much bigger than $\frac{\mathrm{d}N_t}{\mathrm{d}t}$. We can transform Eq.(\ref{num1}) into 
	\begin{equation}\label{num2}
	\left\{
	\begin{array}{l}
		\vspace{2ex}	
		\displaystyle\frac{\mathrm{d}^2 N_t^{\alpha}}{\mathrm{d}t^2}=\displaystyle\Big(-55.1435+1.134632\Phi^{-1}(\alpha)\Big)\frac{\mathrm{d}N_t^{\alpha}}{\mathrm{d}t}\\
		\ \ \ \ \ \ \ \ \ \ \ \ \ \displaystyle+\Big(0.785+2.96798\Phi^{-1}(\alpha)\Big)N_t^{\alpha}\\
		\vspace{2ex}	
		N(0)=1.2157\\	
		N'(0)=0.008.
	\end{array}
	\right.
\end{equation}

When $\alpha>0.4$, $$\Big(-55.1435+1.134632\Phi^{-1}(\alpha)\Big)\frac{\mathrm{d}N_t^{\alpha}}{\mathrm{d}t}+\Big(0.785+2.96798\Phi^{-1}(\alpha)\Big)N_t^{\alpha}$$ is a monotonically increasing function with respect to $N_t^{\alpha}$. Thus, by Theorem \ref{main Formula}, 
for the uncertain nuclear reactor kinetics Eq.\eqref{num}, the inverse uncertainty distribution $\Psi_t^{-1}(\alpha)$ is the solution $N_t^{\alpha}$ of Eq.\eqref{num2}, i.e.,
	\begin{equation*}
		\begin{array}{l}
		\displaystyle\Psi_t^{-1}(\alpha)\!=\!\mathrm{exp}\left\{\frac{(P\!-\!\sqrt{T})}{2}t\right\}\cdot \frac{-0.016\!+\!1.2157 (\sqrt{T}\!+\! P)}{2\sqrt{T}}\\
		\ \ \ \ \ \ \ \ \ \ \ \ \ \  \displaystyle+\mathrm{exp}\left\{\frac{(P\!+\!\sqrt{T})}{2}t\right\}\cdot \frac{0.016\!+\!1.2157(\sqrt{T}\!-\!P) }{2\sqrt{T}}
		\end{array}
	\end{equation*}

where
$$
\begin{aligned}
	&P=-55.1435+1.134632\Phi^{-1}(\alpha),\\
	&R=0.785+2.96798\Phi^{-1}(\alpha),\\
	&T=P^{2}+4R,\\
	&\Phi^{-1}(\alpha)=\frac{\sqrt{3}}{\pi}\mathrm{ln}\frac{\alpha}{1-\alpha},\quad 0.4<\alpha<1.
\end{aligned}
$$

\section{Conclusion}\label{Conclusion}
The comparison theorems of high-order ordinary differential equations are rigorously proved. And the $\alpha$-path of high-order uncertain differential equation is presented. The method to solve a family of high-order uncertain differential equations is proposed including parameter estimation and hypothesis test. Uncertain nuclear reactor kinetics equation is introduced as an example to illustrate this method.


\end{document}